\documentclass[twoside,11pt,reqno]{amsart}
\usepackage{amsmath,amssymb,amscd}

\makeatletter

\@settopoint\oddsidemargin
\@settopoint\marginparwidth
\@settopoint\evensidemargin
\@settopoint\topmargin
\newtheorem{Proposition}{Proposition}[section]
\newtheorem{Lemma}[Proposition]{Lemma}
\newtheorem{Theorem}[Proposition]{Theorem}
\newtheorem{Corollary}[Proposition]{Corollary}

\newtheorem{Remark}[Proposition]{Remark}

\@addtoreset{equation}{section}

\def\k{h}
\def\l{k}
\def\LL{{{\scriptscriptstyle\text{\rm{L}}}}}
\def\F{\mathrm{F}}
\def\L{\mathrm{L}}

\def\C{{\mathbb C}}

\def\0{{\bar 0}}
\def\1{{\bar 1}}

\def\ad{\operatorname{ad}}

\def\gr{\operatorname{gr}}

\def\Mat{M}

\def\sgn{{\operatorname{sgn}}}

\def\underbar{\mathpalette\@underbar}
\def\@underbar#1#2{\settowidth{\@tempdimb}{$#1#2$}\@tempdimb=0.8\@tempdimb
                   \ooalign{$#1#2$\crcr%
                         \hfil\rule[-.5mm]{\@tempdimb}{.4pt}\hfil}}

\newcommand{\bi}{{\text{\boldmath{$i$}}}}
\newcommand{\bj}{{\text{\boldmath{$j$}}}}
\newcommand{\bk}{{\text{\boldmath{$k$}}}}

\def\eps{{\varepsilon}}
\def\phi{{\varphi}}

\begin{document}
\title[Presentations of the Yangian]{\boldmath Parabolic presentations 
of the Yangian $Y(\mathfrak{gl}_n)$}
\author{Jonathan Brundan and Alexander Kleshchev}
\address
{Department of Mathematics\\ University of Oregon\\
Eugene\\ OR~97403, USA}
\email{brundan@darkwing.uoregon.edu, klesh@math.uoregon.edu}
\thanks{Research partially supported by the NSF (grant no. DMS-0139019).}
\thanks{
{\em 2000 Subject Classification}: 17B37.}
\maketitle

\begin{abstract}
We introduce some new presentations for the Yangian 
associated to the Lie algebra $\mathfrak{gl}_n$.
These presentations are
parametrized by tuples of positive integers summing to $n$.
At one extreme, for the tuple $(n)$, the presentation is the usual
RTT presentation of $Y_n$. At the other extreme, for the tuple
$(1^n)$, the presentation is closely related to Drinfeld's presentation.
In general, the presentations are useful for understanding the structure of the
standard parabolic subalgebras of $Y_n$.
\end{abstract}

\section{Introduction}\label{sintro}

Let $Y_n = Y(\mathfrak{gl}_n)$ denote the Yangian associated to the Lie algebra
$\mathfrak{gl}_n$ over the ground field $\C$; see e.g. \cite{D1}, 
\cite[ch.12]{CP} and \cite{MNO}.
In this article, we record some new presentations of $Y_n$
that are adapted to standard parabolic subalgebras.
Let us formulate the main result of the article precisely right away, 
even though the relations appearing in the statement look rather scary.
(Note that in $\S\S$5--6 most of these formulae
are written down in a more concise way in terms of generating series.) Let
$\mathfrak{gl}_\nu = \mathfrak{gl}_{\nu_1} \oplus\cdots\oplus
\mathfrak{gl}_{\nu_m}$ be a standard Levi subalgebra 
of $\mathfrak{gl}_n$,
so $\nu = (\nu_1,\dots,\nu_m)$ is a tuple
of positive integers summing to
 $n$.

\vspace{2mm}
\noindent{\bf Theorem A.}
{\em
The algebra $Y_n$ is generated by elements
$\{D_{a;i,j}^{(r)}, \widetilde{D}_{a;i,j}^{(r)}\}_{1 \leq a \leq m,  1 \leq i,j \leq \nu_a,r \geq 0}$,
$\{E_{a;i,j}^{(r)}\}_{1 \leq a < m, 1 \leq i \leq \nu_a,1 \leq j \leq \nu_{a+1}, r \geq 1}$ and
$\{F_{a;i,j}^{(r)}\}_{1 \leq a < m, 1 \leq i \leq \nu_{a+1},
1 \leq j \leq \nu_a, r \geq 1}$ subject only to the relations
\begin{align}
D_{a;i,j}^{(0)} &= \delta_{i,j},\label{pr1}\\
\sum_{t=0}^r
D_{a;i,p}^{(t)} \widetilde{D}_{a;p,j}^{(r-t)} &= -\delta_{r,0}\delta_{i,j},\label{pr2}\\
[D_{a;i,j}^{(r)}, D_{b;\k,\l}^{(s)}] &=
\delta_{a,b}
\sum_{t=0}^{\min(r,s)-1}
\left(
D_{a;i,\l}^{(r+s-1-t)}D_{a;\k,j}^{(t)} 
-D_{a;i,\l}^{(t)}D_{a;\k,j}^{(r+s-1-t)} \right),\label{pr3}\\
[E_{a;i,j}^{(r)},F_{b;\k,\l}^{(s)}]
&=
\delta_{a,b} \sum_{t=0}^{r+s-1}
\widetilde{D}_{a;i,\l}^{(t)}D_{a+1;\k,j}^{(r+s-1-t)},\label{pr6}\\
[D_{a;i,j}^{(r)}, E_{b;\k,\l}^{(s)}] &=
\delta_{a,b} 
\sum_{t=0}^{r-1} 
D_{a;i,p}^{(t)} E_{a;p,\l}^{(r+s-1-t)}\delta_{\k,j}
- \delta_{a,b+1} \sum_{t=0}^{r-1}
D_{b+1;i,\l}^{(t)} E_{b;\k,j}^{(r+s-1-t)},\label{pr4}\\
[D_{a;i,j}^{(r)}, F_{b;\k,\l}^{(s)}] &=
\delta_{a,b+1} \sum_{t=0}^{r-1}
 F_{b;i,\l}^{(r+s-1-t)}D_{b+1;\k,j}^{(t)}-\delta_{a,b} 
\delta_{i,\l}\sum_{t=0}^{r-1} 
F_{a;\k,p}^{(r+s-1-t)}D_{a;p,j}^{(t)},\label{pr5}
\end{align}\begin{align}
[E_{a;i,j}^{(r)}, E_{a;\k,\l}^{(s)}] 
&=
\sum_{t=1}^{s-1} E_{a;i,\l}^{(t)} E_{a;\k,j}^{(r+s-1-t)}
-\sum_{t=1}^{r-1} E_{a;i,\l}^{(t)} E_{a;\k,j}^{(r+s-1-t)},\label{pr7}\\
[F_{a;i,j}^{(r)}, F_{a;\k,\l}^{(s)}] 
&=
\sum_{t=1}^{r-1} F_{a;i,\l}^{(r+s-1-t)}F_{a;\k,j}^{(t)}-
\sum_{t=1}^{s-1} F_{a;i,\l}^{(r+s-1-t)}F_{a;\k,j}^{(t)},\label{pr8}\\
[E_{a;i,j}^{(r)}, E_{a+1;\k,\l}^{(s+1)}]
&-[E_{a;i,j}^{(r+1)}, E_{a+1;\k,\l}^{(s)}]
=
-E_{a;i,q}^{(r)} E_{a+1;q,\l}^{(s)}\delta_{\k,j},\label{pr9}\\
[F_{a;i,j}^{(r+1)}, F_{a+1;\k,\l}^{(s)}]
&-[F_{a;i,j}^{(r)}, F_{a+1;\k,\l}^{(s+1)}]
=
-\delta_{i,\l}F_{a+1;\k,q}^{(s)}F_{a;q,j}^{(r)},\label{pr10}\\
[E_{a;i,j}^{(r)}, E_{b;\k,\l}^{(s)}] &= 0
\quad\text{ if $b>a+1$ or if $b=a+1$ and $\k \neq j$},\label{pr11}\\
[F_{a;i,j}^{(r)}, F_{b;\k,\l}^{(s)}] &= 0 
\quad\text{ if $b>a+1$ or if $b=a+1$ and $i \neq \l$},\label{pr12}\\
[E_{a;i,j}^{(r)}, [E_{a;\k,\l}^{(s)}, E_{b;f,g}^{(t)}]]
&+
[E_{a;i,j}^{(s)}, [E_{a;\k,\l}^{(r)}, E_{b;f,g}^{(t)}]] = 0
\quad\:\:\text{ if }|a-b|=1,\label{pr13}\\
[F_{a;i,j}^{(r)}, [F_{a;\k,\l}^{(s)}, F_{b;f,g}^{(t)}]]
&+
[F_{a;i,j}^{(s)}, [F_{a;\k,\l}^{(r)}, F_{b;f,g}^{(t)}]] = 0
\,\,\:\:\quad\text{ if }|a-b|=1,\label{pr14}
\end{align}
for all admissible $a,b,f,g,h,i,j,k,r,s,t$.
(By convention the index $p$ resp. $q$ appearing here should be summed over
$1,\dots,\nu_a$ resp. $1,\dots,\nu_{a+1}$.)
}
\vspace{2mm}

In the special case $\nu = (n)$, this presentation is the
RTT presentation of $Y_n$ originating in the work of
Faddeev, Reshetikhin and Takhtadzhyan \cite{FRT}, while if $\nu = (1^n)$ 
the presentation is a variation on Drinfeld's presentation from 
\cite{D3}; see Remark~\ref{drinrem}
for the precise relationship. 
One reason that Drinfeld's presentation is important is because it allows
one to define subalgebras of $Y_n$ which play the role
of the Cartan subalgebra and Borel subalgebra in classical Lie theory.
Our  presentations are well-suited to defining
standard Levi and parabolic subalgebras.
In the notation of Theorem A, 
let $Y_\nu$, $Y_\nu^+$ and $Y_\nu^-$ denote the subalgebras
of $Y_n$ generated by all the $D_{a;i,j}^{(r)}$'s,
the $E_{a;i,j}^{(r)}$'s or the $F_{a;i,j}^{(r)}$'s, respectively. 
Then, $Y_\nu = Y(\mathfrak{gl}_\nu)$ is the {\em standard Levi subalgebra} 
of $Y_n$, isomorphic to
$Y_{\nu_1} \otimes \cdots \otimes Y_{\nu_m}$.
The {\em standard parabolic subalgebras} $Y_\nu^\sharp$
and $Y_\nu^\flat$
of $Y_n$ are the
subalgebras generated by $Y_\nu, Y_\nu^+$ and by
$Y_\nu,Y_\nu^-$ respectively; see also Remark~\ref{new} where these subalgebras
are defined directly in terms of the Drinfeld presentation.

The next theorem
gives a PBW basis for each of these algebras.
To write it down,
define elements
$E_{a,b;i,j}^{(r)}$ and $F_{a,b;j,i}^{(r)}$ for $1 \leq a < b \leq m$
and $1 \leq i \leq \nu_a, 1 \leq j \leq \nu_b$ 
by setting $E_{a,a+1;i,j}^{(r)} := E_{a;i,j}^{(r)},
F_{a,a+1;j,i}^{(r)} := F_{a;j,i}^{(r)}$ and then 
inductively defining 
$$
E_{a,b;i,j}^{(r)} := [E_{a,b-1;i,k}^{(r)}, E_{b-1;k,j}^{(1)}],
\qquad
F_{a,b;j,i}^{(r)} := [F_{b-1;j,k}^{(1)}, F_{a,b-1;k,i}^{(r)}],
$$
if $b > a+1$, where $1 \leq k \leq \nu_{b-1}$.
The relations imply that these definitions
are independent of the particular choice of $k$; see (\ref{indofk}).

\vspace{2mm}

\noindent
{\bf Theorem B.}
{\em {\rm(i)} The set of all monomials in 
$\{D_{a;i,j}^{(r)}\}_{a=1,\dots,m, 1 \leq i,j \leq \nu_a,r\geq 1}$ 
taken in some
\begin{itemize}
\item[]
fixed order forms a basis for $Y_\nu$.
\item[(ii)] The set of all monomials in 
$\{E_{a,b;i,j}^{(r)}\}_{1 \leq a < b \leq m, 1 \leq i \leq \nu_a, 1 \leq j \leq \nu_b, r \geq 1}$ taken in some fixed
order forms a basis for $Y_\nu^+$.
\item[(iii)] The set of all monomials in 
$\{F_{a,b;i,j}^{(r)}\}_{1 \leq a < b \leq m, 
1 \leq i \leq \nu_b, 1 \leq j \leq \nu_a, r \geq 1}$ taken in some fixed
order forms a basis for $Y_\nu^-$.
\item[(iv)] The set of all monomials in the union of
the elements listed in (i),(ii) and (iii)
taken in some fixed order forms a basis for $Y_n$.
\end{itemize}}

\vspace{2mm}

This theorem implies in particular that 
the natural multiplication maps $Y_\nu^- \otimes Y_\nu \otimes Y_\nu^+ 
\rightarrow Y_n$, 
$Y_\nu \otimes Y_\nu^+ \rightarrow Y_\nu^\sharp$ and
$Y_\nu^- \otimes Y_\nu \rightarrow Y_\nu^\flat$
are vector space isomorphisms. Moreover, 
there are natural
projections $Y_\nu^{\sharp} \twoheadrightarrow Y_\nu$ and $Y_\nu^\flat 
\twoheadrightarrow
Y_\nu$ with kernels generated
by all $E_{a;i,j}^{(r)}$ and by all 
$F_{a;i,j}^{(r)}$ respectively.

The rest of the article is organized as follows.
To start with,
$\S\S$\ref{srtt}-\ref{spbw} are expository, giving the necessary
definitions and a proof of the PBW theorem for $Y_n$.
In $\S$\ref{slevi}, we define Levi subalgebras.
Then in $\S$\ref{sdrinfeld} we give a complete proof of the
equivalence of the RTT presentation of $Y_n$ with Drinfeld's presentation,
filling a gap in the literature.
Note though that a proof of the analogous but harder 
result in the quantum affine case can be found in work of Ding and Frenkel 
\cite{DF}, and the basic trick of considering certain Gauss factorizations is the same here.
The main theorems are proved in
$\S$\ref{sparabolic}. The argument involves {\em partial}
Gauss factorizations, an idea already exploited by Ding \cite{Ding} 
to study the embedding of $U_q(\widehat{\mathfrak{gl}}_{n-1})$
in $U_q(\widehat{\mathfrak{gl}}_n)$
in the quantum affine setting.
The remaining two sections
are again expository in nature:
in $\S$\ref{scent} we record proofs of some known results
about centers and centralizers and in $\S$\ref{squantum} we explain
the relationship between
our approach and the {\em quantum determinants} which are used
to define the Drinfeld generators elsewhere in the literature.

The results of this article play a central role in 
\cite{BK}, where we derive
generators and relations for the finite $W$-algebras
associated to nilpotent matrices in the general linear Lie algebras.

\vspace{1mm}
\noindent
{\em Acknowledgements.}
The second author would like to thank Arun Ram for stimulating 
conversations.

\vspace{1mm}
\noindent
{\em Notation.}
Throughout the article, we work over the ground field $\C$.
We write $M_n$ for the associative algebra of all $n \times n$ matrices over $\C$, and $\mathfrak{gl}_n$ for the corresponding Lie algebra.
The $ij$-matrix unit is denoted $e_{i,j}$.

\section{RTT presentation}\label{srtt}
To define the {\em Yangian} $Y_n = Y(\mathfrak{gl}_n)$
we use the RTT formalism; see \cite[ch. 11]{ES} or \cite{FRT}.
Our basic reference for this material in the case of the Yangian
is \cite[$\S$1]{MNO}.
Let
\begin{equation}
R(u) = u - \sum_{i,j=1}^n e_{i,j} \otimes e_{j,i} \in \Mat_n \otimes 
\Mat_n[u]
\end{equation}
denote Yang's $R$-matrix with parameter $u$. 
This satisfies the QYBE with spectral parameters:
\begin{equation}
R^{[1,2]}(u-v) R^{[1,3]}(u) R^{[2,3]}(v) =
R^{[2,3]}(v) R^{[1,3]}(u) R^{[1,2]}(u-v),
\end{equation}
equality in $\Mat_n \otimes \Mat_n \otimes \Mat_n[u,v]$. 
The superscripts in square brackets here and later on
indicate the embedding of a smaller tensor into a bigger tensor,
inserting the identity into all other tensor positions.
Now, $Y_n$ is defined to be the associative algebra on generators
$\{T_{i,j}^{(r)}\}_{1 \leq i,j \leq n, r \geq 1}$ subject to certain
relations.
To write down these relations, let
$$
T_{i,j}(u) := \sum_{r \geq 0} T_{i,j}^{(r)} u^{-r} \in Y_n [[u^{-1}]]
$$
where $T_{i,j}^{(0)} := \delta_{i,j}$, and
\begin{equation}
T(u) := \sum_{i,j=1}^n e_{i,j} \otimes T_{i,j}(u)
\in \Mat_n \otimes Y_n[[u^{-1}]].
\end{equation}
We often think of $T(u)$ as an $n \times n$ matrix with $ij$-entry
$T_{i,j}(u)$.
Now the relations are given by the equation
\begin{equation}\label{rmdef}
R^{[1,2]}(u-v) T^{[1,3]}(u) T^{[2,3]}(v) =
 T^{[2,3]}(v) T^{[1,3]}(u) R^{[1,2]}(u-v),
\end{equation}
equality in $\Mat_n \otimes \Mat_n \otimes Y_n((u^{-1},v^{-1}))$
(the localization of $\Mat_n \otimes \Mat_n \otimes Y_n[[u^{-1}, v^{-1}]]$
at the multiplicative set consisting of the non-zero elements of
$\C[[u^{-1} ,v^{-1}]]$).
Equating $e_{i,j} \otimes e_{\k,\l} \otimes ?$-coefficients on either side of
(\ref{rmdef}), the relations are equivalent to
\begin{equation}\label{trel0}
(u-v)[T_{i,j}(u), T_{\k,\l}(v)] = T_{\k,j}(u)T_{i,\l}(v)  - T_{\k,j}(v)T_{i,\l}(u) .
\end{equation}
Swapping $i$ with $\k$, $j$ with $\l$ and $u$ with $v$, we get equivalently
that
\begin{equation}\label{trel}
(u-v)[T_{i,j}(u), T_{\k,\l}(v)] = T_{i,\l}(v) T_{\k,j}(u) - T_{i,\l}(u) T_{\k,j}(v).
\end{equation}
Yet another formulation of the relations is given by 
\begin{equation}\label{mr}
[T_{i,j}^{(r)}, T_{\k,\l}^{(s)}]
= \sum_{t=0}^{\min(r,s)-1}
\left( 
T_{i,\l}^{(r+s-1-t)}T_{\k,j}^{(t)} -
T_{i,\l}^{(t)}T_{\k,j}^{(r+s-1-t)}
\right)
\end{equation}
for every $1 \leq h,i,j,k \leq n$ and $r,s \geq 0$; see 
\cite[Proposition 1.2]{MNO}.

Using (\ref{rmdef}), one checks
that the following are (anti)automorphisms
of $Y_n$;
see \cite[Proposition 1.12]{MNO}.
\begin{itemize}
\item[(A1)] ({\em Translation}) For $c \in \C$, let $\eta_c:Y_n 
\rightarrow Y_n$
be the automorphism defined from $\eta_c^{[2]}(T(u)) = T(u+c)$, 
i.e. 
$\eta_c(T_{i,j}^{(r)}) = \sum_{s=0}^{r-1}
\binom{r-1}{s} (-c)^s T_{i,j}^{(r-s)}$.
\item[(A2)] ({\em Multiplication by a
power series}) For $f(u) \in 1 + u^{-1} \C[[u^{-1}]]$, let 
$\mu_f:Y_n \rightarrow Y_n$ be the automorphism defined from
$\mu_f^{[2]}(T(u)) = f(u) T(u)$, i.e. 
$\mu_f(T_{i,j}^{(r)}) = \sum_{s=0}^r a_s T_{i,j}^{(r-s)}$
if $f(u) = \sum_{s \geq 0} a_s u^{-s}$.
\item[(A3)] ({\em Sign change}) 
Let $\sigma:Y_n\rightarrow Y_n$ 
be the antiautomorphism of order 2
defined from
$\sigma^{[2]}(T(u)) = T(-u)$, i.e.
$\sigma(T_{i,j}^{(r)}) = (-1)^r T_{i,j}^{(r)}$.
\item[(A4)] ({\em Transposition})
Let $\tau:Y_n\rightarrow Y_n$
be the antiautomorphism of order 2
defined from
$\tau^{[2]}(T(u)) = (T(u))^{t}$ (transpose matrix), i.e. 
$\tau(T_{i,j}^{(r)}) = T_{j,i}^{(r)}$.
\item[(A5)] ({\em Inversion})
Let $\omega:Y_n \rightarrow Y_n$ be the
automorphism of order $2$ defined from the equation
$\omega^{[2]}(T(u)) = T(-u)^{-1}$.
\end{itemize}

The Yangian $Y_n$ is a Hopf algebra with
comultiplication $\Delta:Y_n \rightarrow Y_n \otimes Y_n$,
counit $\eps:Y_n\rightarrow\C$ and antipode $S:Y_n\rightarrow Y_n$
defined from
$$
\Delta^{[2]}(T(u)) = T^{[1,2]}(u) T^{[1,3]}(u),\quad
\eps^{[2]}(T(u)) = I,\quad S^{[2]}(T(u)) = T(u)^{-1},
$$
equalities written in $\Mat_n\otimes Y_n\otimes Y_n[[u^{-1}]],
\Mat_n[[u^{-1}]]$ and $\Mat_n\otimes Y_n[[u^{-1}]]$ respectively.
Note that $S = \omega \circ \sigma$. 
The involutions $\omega$ and $\sigma$ do {\em not} commute,
so $S$ is not of order 2; a precise description of $S^2$ is given in 
\cite[Theorem 5.11]{MNO} or
Corollary~\ref{s2} below.
Since it arises quite often, we let
$\widetilde{T}_{i,j}(u) = \sum_{r \geq 0} \widetilde T_{i,j}^{(r)} u^{-r}
:= -S(T_{i,j}(u))$,
so
\begin{align}
\widetilde{T}(u) := \sum_{i,j=1}^n e_{i,j} \otimes \widetilde{T}_{i,j}(u)  = 
-T(u)^{-1}.
\end{align}
To work out commutation relations between $T_{i,j}(u)$
and $\widetilde T_{\k,\l}(v)$, it is useful to rewrite the RTT presentation in the
form
\begin{equation}\label{rmdefalt}
\widetilde T^{[2,3]}(v) R^{[1,2]}(u-v) T^{[1,3]}(u)  =
 T^{[1,3]}(u) R^{[1,2]}(u-v)\widetilde T^{[2,3]}(v) .
\end{equation}
We record \cite[Lemma 1.1]{NT}:

\begin{Lemma}\label{cim}
Given $i \neq \l$ and $\k \neq j$, $T_{i,j}^{(r)}$ and
$\widetilde T_{\k,\l}^{(s)}$ commute for all $r,s \geq 1$.
\end{Lemma}

\begin{proof}
Compute the $e_{i,j} \otimes e_{\k,\l}$-coefficients on each side
of (\ref{rmdefalt}).
\end{proof}

We often work with the {\em canonical filtration}
\begin{equation}\label{filt1}
\F_0 Y_n \subseteq \F_1 Y_n \subseteq 
\F_2 Y_n \subseteq \cdots
\end{equation}
on $Y_n$
defined by declaring that the generator $T_{i,j}^{(r)}$ is of degree $r$
for each $r \geq 1$, i.e.
$\F_d Y_n$ is the span of all monomials of the form
$T_{i_1,j_1}^{(r_1)} \cdots T_{i_s,j_s}^{(r_s)}$
with {total degree} $r_1+\cdots+r_s$ at most $d$.
It is easy to see using (\ref{mr}) 
that the associated graded algebra
$\operatorname{gr} Y_n$ is commutative.
From this one deduces by induction on $d$ that $\F_d Y_n$ 
is already spanned by the set of all monomials 
of total degree at most $d$
in the elements 
$\{T_{i,j}^{(r)}\}_{1\leq i,j \leq n, r \geq 1}$ 
{\em taken in some fixed order}.
In fact it is known that such ordered monomials are linearly independent,
hence the set of all monomials in the 
elements
$\{T_{i,j}^{(r)}\}_{1\leq i,j \leq n, r \geq 1}$ taken in some fixed order
gives a basis for $Y_n$;
see \cite[Corollary 1.23]{MNO} or Corollary~\ref{pbwcor} below.
In other words, the associated graded algebra $\gr Y_n$
is the free commutative algebra on generators
$\{\gr_r T_{i,j}^{(r)}\}_{1 \leq i,j \leq n, r \geq 1}$.

There is a second important filtration which we call the
{\em loop filtration}
\begin{equation}\label{filt2}
\L_0 Y_n \subseteq \L_1 Y_n \subseteq \L_2 Y_n \subseteq \cdots
\end{equation}
defined by declaring that the generator
$T_{i,j}^{(r)}$ is of degree $(r-1)$
for each $r \geq 1$.
We denote the associated graded algebra by $\gr^\LL Y_n$.
Let $\mathfrak{gl}_n[t]$ denote the Lie algebra
$\mathfrak{gl}_n \otimes \C[t]$ with basis
$\{e_{i,j}t^r\}_{1 \leq i,j \leq n, r \geq 0}$.
By the relations (\ref{mr}), there is a surjective homomorphism
$U(\mathfrak{gl}_n[t]) \twoheadrightarrow \gr^\LL Y_n$
mapping $e_{i,j} t^r$ to $\gr^\LL_{r} T_{i,j}^{(r+1)}$
for each $1 \leq i,j  \leq n, r \geq 0$.
By the PBW theorem for $Y_n$ described in the previous paragraph, this 
map is actually an isomorphism, hence
$\gr^\LL Y_n \cong U(\mathfrak{gl}_n[t])$; see also \cite[Theorem 1.26]{MNO}
where this argument is explained in more detail.

\section{PBW theorem}\label{spbw}

In this section, we give a proof of the PBW theorem for $Y_n$
different from the one in \cite{MNO}.
It was inspired by the realization of the Yangian found in
\cite{BR}.
Let $U(\mathfrak{gl}_n)$ denote the universal enveloping algebra
of the Lie algebra $\mathfrak{gl}_n$.
We have the {\em evaluation homomorphism} 
\begin{equation}\label{ev0}
\kappa_1:Y_n \rightarrow U(\mathfrak{gl}_n),\quad
T_{i,j}^{(1)} \mapsto 
\left\{
\begin{array}{ll}
e_{i,j}&\hbox{if $r = 1$,}\\
0&\hbox{if $r > 1$.}
\end{array}\right.
\end{equation}
More generally for $l \geq 1$, consider the homomorphism
\begin{equation}
\kappa_{l} := 
\kappa_1 \otimes\cdots\otimes \kappa_1 \circ \Delta^{(l)}:Y_n
\rightarrow U(\mathfrak{gl}_n)^{\otimes l},
\end{equation}
where $\Delta^{(l)}:Y_n \rightarrow Y_n^{\otimes l}$ 
denotes the $l$th iterated comultiplication.
We define the algebra
$Y_{n,l}$ to be the image $\kappa_{l}(Y_n)$ of $Y_n$
under this homomorphism.
Writing $e_{i,j}^{[s]}$ for the
element $1^{\otimes{(s-1)}} \otimes e_{i,j} \otimes 1^{\otimes(l-s)}
\in U(\mathfrak{gl}_n)^{\otimes l}$, we have by the definition of
$\Delta$ that
\begin{equation}
\kappa_{l}(T_{i,j}^{(r)})
=
\sum_{1 \leq s_1 < \cdots < s_r \leq l}
\sum_{\substack{1 \leq i_0,\cdots,i_r \leq n \\ i_0 = i, i_r = j}} 
e_{i_0, i_1}^{[s_1]}
e_{i_1, i_2}^{[s_2]}
\cdots e_{i_{r-1}, i_r}^{[s_r]}.
\end{equation}
In particular, we see from this that $\kappa_{l}(T_{i,j}^{(r)}) = 0$ for all
$r > l$.

\begin{Theorem}\label{truncthm}
The set of all 
monomials in the elements $\{\kappa_{l}(T_{i,j}^{(r)})\}_{
1 \leq i,j \leq n, r =1,\dots,l}$ taken in some fixed order
forms a basis for $Y_{n,l}$.
\end{Theorem}

\begin{proof}
It is obvious that such monomials span $Y_{n,l}$, so we just
have to show that they are linearly independent.
Consider the standard filtration
$$
\F_0 U(\mathfrak{gl}_n)^{\otimes l}
\subseteq \F_1 U(\mathfrak{gl}_n)^{\otimes l} \subseteq \F_2 U(\mathfrak{gl}_n)^{\otimes l} \subseteq \cdots
$$
on $U(\mathfrak{gl}_n)^{\otimes l}$, 
so each generator $e_{i,j}^{[r]}$ is of degree $1$. 
The associated graded algebra $\gr U(\mathfrak{gl}_n)^{\otimes l}$
is the free polynomial algebra
on generators $x_{i,j}^{[r]} := \gr_1 e_{i,j}^{[r]}$. 
Let $y_{i,j}^{(r)} := \gr_r \kappa_{l}(T_{i,j}^{(r)})$,
i.e.
$$
y_{i,j}^{(r)}
=
\sum_{1 \leq s_1 < \cdots < s_r \leq l}
\sum_{\substack{1 \leq i_0,\cdots,i_r \leq n \\ i_0 = i, i_r = j}} 
x_{i_0, i_1}^{[s_1]}x_{i_1, i_2}^{[s_2]}\cdots x_{i_{r-1}, i_r}^{[s_r]}.
$$
To complete the proof of the theorem, we show that
the elements $\{y_{i,j}^{(r)}\}_{1 \leq i,j \leq n, r =1,\dots,l}$
are algebraically independent.

Let us identify $\gr U(\mathfrak{gl}_n)^{\otimes l}$
with the coordinate algebra 
$\C[M_n^{\times l}]$ of the affine variety $M_n^{\times l}$
of $l$-tuples $(A_1,\dots, A_l)$ of
$n \times n$ matrices, so that
$x_{i,j}^{[r]}$ is the function picking out the $ij$-entry of the
$r$th matrix $A_r$.
Let $\theta:M_n^{\times l} \rightarrow M_n^{\times l}$
be the morphism defined by 
$(A_1,\dots,A_l) \mapsto (B_1,\dots,B_l)$, where
$B_r$
is the $r$th elementary symmetric function 
$$
e_r(A_1,\dots,A_l)
:=
\sum_{1 \leq s_1 < \cdots < s_r \leq l} A_{s_1} \cdots A_{s_r}
$$ 
in the matrices $A_1,\dots,A_l$.
The comorphism $\theta^*$ maps $x_{i,j}^{[r]}$
to $y_{i,j}^{(r)}$. 
So to show that the $y_{i,j}^{(r)}$ are algebraically
independent, we need to show that $\theta^*$ is injective,
i.e. that $\theta$ is a dominant morphism of affine varieties.
For this it suffices to show that the differential of $\theta$
is surjective at some point $x \in M_n^{\times l}$.

Pick pairwise distinct scalars $c_1,\dots,c_l \in \C$ and consider 
$x := (c_1 I_n, \dots, c_l I_n)$.
Identifying the tangent space $T_x(M_n^{\times l})$
with the vector space $M_n^{\oplus l}$,
a calculation shows that the differential $d \theta_x$ maps 
$(A_1,\dots,A_l)$ to $(B_1,\dots,B_l)$ where 
$$
B_r = \sum_{s=1}^l e_{r-1}(c_1,\dots,\widehat{c_s},\dots,c_l) A_s.
$$
Here, $e_{r-1}(c_1,\dots,\widehat{c_s},\dots,c_l)$ 
denotes the $(r-1)$th elementary symmetric function
in the scalars $c_1,\dots,c_l$ excluding $c_s$.
We just need to show this linear map is surjective, for which it
clearly suffices to consider the case $n=1$. But in that case
its determinant is the Vandermonde determinant
$\prod_{1 \leq r < s \leq l} (c_s-c_r)$, so it is non-zero by the
choice of the scalars $c_1,\dots,c_l$.
\end{proof}

\begin{Corollary}\label{pbwcor}
The 
set of all 
monomials in the elements $\{T_{i,j}^{(r)}\}_{
1 \leq i,j \leq n, r \geq 1}$ taken in some fixed order
forms a basis for $Y_{n}$.
\end{Corollary}

\begin{proof}
We have already observed that such monomials span $Y_n$.
The fact that they are linearly independent follows from
Theorem~\ref{truncthm} by taking sufficiently large $l$.
\end{proof}

\begin{Corollary}\label{trunccor}
The kernel of $\kappa_{l}:Y_n \twoheadrightarrow Y_{n,l}$
is the two-sided ideal of $Y_n$
generated by the elements
$\{T_{i,j}^{(r)}\}_{1 \leq i,j \leq n, r > l}$.
\end{Corollary}

\begin{proof} 
Let $I$ denote the two-sided ideal of $Y_n$
generated by $\{T_{i,j}^{(r)}\}_{1 \leq i,j \leq n, r > l}$.
It is obvious that $\kappa_{l}$ induces a map  
$\bar\kappa_{l}: Y_n / I \twoheadrightarrow Y_{n,l}$.
Since $Y_n / I$ is spanned by the
set of all monomials in the elements $\{T_{i,j}^{(r)}+I\}_{
1 \leq i,j \leq n, r =1,\dots,l}$ taken in some fixed order,
Theorem~\ref{truncthm} now implies that $\bar\kappa_{l}$ is an isomorphism.
\end{proof}

The first of these corollaries proves the PBW theorem for $Y_n$.
The second corollary shows that the algebra $Y_{n,l}$ is 
the {\em Yangian of level $l$} introduced 
by Cherednik \cite{Ch0,Ch}. Moreover, by Corollary~\ref{trunccor},
the maps $\kappa_{l}$ induce an inverse system
$$
Y_{n,1} \twoheadleftarrow Y_{n,2} \twoheadleftarrow \cdots
$$
of filtered algebras, where each $Y_{n,l}$ is filtered 
by the canonical filtration
defined by declaring
that the generators $\kappa_l (T_{i,j}^{(r)})$ are of degree $r$.
It is easy to see using Theorem~\ref{truncthm} and Corollary~\ref{pbwcor}
that the Yangian $Y_n$ is the 
inverse limit $\varprojlim Y_{n,l}$ of this system 
taken in the
category of filtered algebras. This gives a concrete realization of
the Yangian.

\section{Levi subalgebras}\label{slevi}
Our exposition is biased towards 
the {\em standard embedding} $Y_{n} \hookrightarrow Y_{n+1}$
under which $T_{i,j}^{(r)} \in Y_{n}$ maps to the element with the
same name in $Y_{n+1}$. 
We warn the reader that the element
$\widetilde{T}_{i,j}^{(r)} \in Y_n$ does {\em not} map to the element with
the same name in $Y_{n+1}$!
The standard embeddings define a tower of algebras
\begin{equation}\label{tower}
Y_1 \subset Y_2 \subset Y_3 \subset \cdots
\end{equation}
which will be implicit in our work from now on.
Since most of the automorphisms of the Yangian defined in 
$\S$\ref{srtt} do {\em not} commute with the standard embeddings,
we sometimes add a subscript to clarify notation;
for example we write $\omega_n:Y_n \rightarrow Y_n$
for the automorphism $\omega$ if confusion seems likely.

For $m \geq 0$, we let 
$\phi_m:Y_n \hookrightarrow Y_{m+n}$ denote the obvious injective
algebra homomorphism mapping $T_{i,j}^{(r)} \in Y_n$ to $T_{m+i,m+j}^{(r)} 
\in Y_{m+n}$
for each $1 \leq i,j \leq n, r \geq 1$.
Following \cite{NT}, we define another injective algebra homomorphism
$\psi_m:Y_n \hookrightarrow Y_{m+n}$ by
\begin{equation}\label{pdidef}
\psi_{m} := \omega_{m+n} \circ \phi_m \circ \omega_n.
\end{equation}
Observe that $\psi_m$ maps $\widetilde{T}_{i,j}^{(r)} \in Y_n$
to $\widetilde{T}_{m+i,m+j}^{(r)} \in Y_{m+n}$. So
the subalgebra $\psi_m(Y_n)$ of $Y_{m+n}$ is generated by the elements
$\{\widetilde{T}_{m+i,m+j}^{(r)}\}_{1 \leq i,j \leq n, r \geq 1}$.
Given this, the following lemma is 
an immediate consequence of Lemma~\ref{cim}.

\begin{Lemma}\label{cent}
The subalgebras $Y_m$ and $\psi_m(Y_n)$ of $Y_{m+n}$ centralize each other.
\end{Lemma}

Let us give another description of the map $\psi_m$ 
in terms of the quasi-determinants of Gelfand and Retakh; 
see e.g. \cite[$\S$2.2]{GKLLRT}.
Suppose that $A, B, C$ and $D$ are $m \times m$, $m \times n$, $n \times m$
and $n \times n$ matrices respectively with entries in some ring $R$.
Assuming that the matrix $A$ is invertible, we define
\begin{equation}\label{qdetdef}
\left|
\begin{array}{cc}
A&B\\
C&
\hbox{\begin{tabular}{|c|}\hline$D$\\\hline\end{tabular}}
\end{array}
\right| := D - C A^{-1} B.
\end{equation}
Then:

\begin{Lemma}\label{qdet}
For any $1 \leq i,j \leq n$, 
$$
\psi_m(T_{i,j}(u))
=
\left|
\begin{array}{cccc}
T_{1,1}(u) & \cdots &T_{1,m}(u)&T_{1,m+j}(u)\\
\vdots & \ddots &\vdots&\vdots\\
T_{m,1}(u) & \cdots & T_{m,m}(u)&T_{m,m+j}(u)\\
T_{m+i,1}(u) & \cdots & T_{m+i,m}(u)&
\hbox{\begin{tabular}{|c|}\hline$T_{m+i,m+j}(u)$\\\hline\end{tabular}}
\end{array}
\right|.
$$
\end{Lemma}

\begin{proof}
Let $T(u)$ denote the matrix $\left( T_{i,j}(u) \right)_{1 \leq i,j \leq n}$
with entries in $Y_n[[u^{-1}]]$ as usual
and let $\widetilde{T}(u) := 
-T(u)^{-1}$.
Also define the matrices
\begin{align*}
A(u) &= \left( T_{i,j}(u) \right)_{1 \leq i,j \leq m},
&B(u) &= \left( T_{i,j}(u) \right)_{1 \leq i \leq m, m+1 \leq j \leq m+n},\\
C(u) &= \left( T_{i,j}(u) \right)_{m+1 \leq i \leq m+n, 1 \leq j \leq m},
&D(u) &= \left( T_{i,j}(u) \right)_{m+1 \leq i,j \leq m+n}
\end{align*}
with entries in $Y_{m+n}[[u^{-1}]]$ and let
$$
\left(
\begin{array}{cc}
\widetilde{A}(u)&\widetilde{B}(u)\\
\widetilde{C}(u)&\widetilde{D}(u)\\
\end{array}
\right)
:=
-\left(
\begin{array}{cc}
A(u)&B(u)\\
C(u)&D(u)\\
\end{array}
\right)^{-1}.
$$
By block multiplication, one checks the classical identities
\begin{align*}
\widetilde{A}(u) &= -\left(A(u)-B(u)D(u)^{-1}C(u)\right)^{-1},\\
\widetilde{B}(u) &=  A(u)^{-1} B(u) \left(D(u)-C(u)A(u)^{-1}B(u)\right)^{-1},\\
\widetilde{C}(u) &=  D(u)^{-1} C(u) \left(A(u)-B(u)D(u)^{-1} C(u)\right)^{-1},\\
\widetilde{D}(u) &= -\left(D(u) - C(u) A(u)^{-1} B(u)\right)^{-1}.
\end{align*}
Now, by its definition, the homomorphism $\psi_m$ maps $\widetilde{T}(u)$ 
to $\widetilde{D}(u)$.
Hence, it maps $T(u)$ to $-\widetilde D(u)^{-1}
=D(u) - C(u) A(u)^{-1} B(u)$.
The lemma follows from this on computing $ij$-entries.
\end{proof}

The description of $\psi_m(T_{i,j}(u))$ given by Lemma~\ref{qdet}
does not depend on $n$. This means that the maps $\psi_m$ are compatible
with the standard embeddings, in the sense that the following diagram commutes
\begin{equation}\label{unam}
\begin{CD}
Y_1 &@>>>&Y_2 &@>>> &Y_{3} &@>>>&\cdots\\
@V\psi_m VV&&@V\psi_m VV&&@V\psi_m VV\\
Y_{m+1}&@>>>&Y_{m+2} &@>>> &Y_{m+3}&@>>>&\cdots
\end{CD}
\end{equation}
where the horizontal maps are standard embeddings. So our notation for the 
map $\psi_m$ is unambiguous as $n$ varies.
We also note that
\begin{equation}\label{comp}
\psi_m \circ \psi_{m'} = \psi_{m+m'}
\end{equation}
for any $m,m' \geq 0$, which is an obvious consequence of 
our original definition.

Now we can define the {\em standard Levi subalgebras}
of $Y_n$. Given a tuple $\nu = (\nu_1,\dots,\nu_m)$ of positive integers
summing to $n$, define $Y_\nu$ to be the subalgebra
\begin{equation}
Y_\nu := Y_{\nu_1} \psi_{\nu_1}(Y_{\nu_2}) \psi_{\nu_1+\nu_2} (Y_{\nu_3}) \cdots
\psi_{\nu_1+\cdots+\nu_{m-1}} (Y_{\nu_m})
\end{equation}
of $Y_n$. For $a=1,\dots,m$ and $1 \leq i,j \leq \nu_a$, we let
\begin{equation}
D_{a;i,j}(u) = \sum_{r \geq 0} D_{a;i,j}^{(r)} u^{-r} :=
\psi_{\nu_1+\cdots+\nu_{a-1}} (T_{i,j}(u)).
\end{equation}
By Lemma~\ref{cent} and induction on $m$, 
the various ``blocks'' of $Y_\nu$ centralize
each other, hence the map
$$
\psi_0 \bar \otimes \psi_{\nu_1} \bar \otimes\cdots\bar \otimes 
\psi_{\nu_1+\cdots+\nu_{m-1}}:
Y_{\nu_1} \otimes Y_{\nu_2} \otimes\cdots\otimes Y_{\nu_m} \rightarrow Y_\nu
$$
is an algebra isomorphism.
This means that 
the elements $\{D_{a;i,j}^{(r)}\}_{1 \leq a \leq m, 1 \leq i,j \leq \nu_a, r \geq 1}$ generate $Y_\nu$ subject only to the following relations:
\begin{equation}\label{levirel}
[D_{a;i,j}^{(r)}, D_{b;\k,\l}^{(s)}]
= \delta_{a,b} 
\sum_{t=0}^{\min(r,s)-1}
\left(D_{a;i,\l}^{(r+s-1-t)}D_{a;\k,j}^{(t)} 
-
D_{a;i,\l}^{(t)}D_{a;\k,j}^{(r+s-1-t)} \right)
\end{equation}
where $D_{a;i,j}^{(0)} := \delta_{i,j}$.
The special case $\nu = (1^n)$ is particularly important:
$Y_{(1^n)} \cong Y_1 \otimes \cdots\otimes Y_1$
is a commutative subalgebra of $Y_n$ which plays the role
of Cartan subalgebra.

\section{Drinfeld presentation}\label{sdrinfeld}

In this section, we introduce the Drinfeld generators
of $Y_n$ following the approach of \cite{DF} and \cite[Appendix B]{Iohara}.
Since the leading minors of the matrix $T(u)$ are invertible, 
it possesses a Gauss factorization
\begin{equation}\label{gfact}
T(u) = F(u) D(u) E(u)
\end{equation}
for unique matrices
$$
D(u) = \left(
\begin{array}{cccc}
D_{1}(u) & 0&\cdots&0\\
0 & D_{2}(u) &\cdots&0\\
\vdots&\vdots&\ddots&\vdots\\
0&0 &\cdots&D_{n}(u)
\end{array}
\right),
$$$$
E(u) = 
\left(
\begin{array}{cccc}
1 & E_{1,2}(u) &\cdots&E_{1,n}(u)\\
0 & 1 &\cdots&E_{2,n}(u)\\
\vdots&\vdots&\ddots&\vdots\\
0&0 &\cdots&1
\end{array}
\right),\:
F(u) = \left(
\begin{array}{cccc}
1 & 0 &\cdots&0\\
F_{1,2}(u) & 1 &\cdots&0\\
\vdots&\vdots&\ddots&\vdots\\
F_{1,n}(u)&F_{2,n}(u) &\cdots&1
\end{array}
\right).
$$
This defines power series
$D_i(u) = \sum_{r \geq 0} D_i^{(r)} u^{-r}$, 
$E_{i,j}(u) = \sum_{r \geq 1} E_{i,j}^{(r)} u^{-r}$ and 
$F_{i,j}(u) = \sum_{r \geq 1} F_{i,j}^{(r)} u^{-r}$.
Let $E_i(u) = \sum_{r \geq 1} E_i^{(r)} u^{-r}
:= E_{i,i+1}(u)$ and $F_i(u) = \sum_{r \geq 1} F_i^{(r)} u^{-r}:= 
F_{i,i+1}(u)$ for short.
Also let $\widetilde D_i(u) =\sum_{r \geq 0} \widetilde{D}_i^{(r)} u^{-r}
:= -D_i(u)^{-1}$. 

In terms of quasi-determinants, we have the following more explicit
descriptions; see \cite[Theorem 4.4]{GR} or
\cite[Theorem 2.2.6]{GR2}.
\begin{equation}\label{qd1}
\!\!\!\!\!\!\!\!\!\!\!\!\!\!\!\!\!\!\!\!
D_i(u) = 
\left|
\begin{array}{cccc}
T_{1,1}(u) & \cdots & T_{1,i-1}(u)&T_{1,i}(u)\\
\vdots & \ddots &\vdots&\vdots\\
T_{i-1,1}(u)&\cdots&T_{i-1,i-1}(u)&T_{i-1,i}(u)\\
T_{i,1}(u) & \cdots & T_{i,i-1}(u)&
\hbox{\begin{tabular}{|c|}\hline$T_{i,i}(u)$\\\hline\end{tabular}}
\end{array}
\right|,
\end{equation}
\begin{equation}
\label{qd2}
E_{i,j}(u) = 
D_i(u)^{-1} \left|
\begin{array}{cccc}
T_{1,1}(u) & \cdots &T_{1,i-1}(u)& T_{1,j}(u)\\
\vdots & \ddots &\vdots&\vdots\\
T_{i-1,1}(u) & \cdots & T_{i-1,i-1}(u)&T_{i-1,j}(u)\\
T_{i,1}(u) & \cdots & T_{i,i-1}(u)&
\hbox{\begin{tabular}{|c|}\hline$T_{i,j}(u)$\\\hline\end{tabular}}
\end{array}
\right|,
\end{equation}
\begin{equation}\label{qd3}
F_{i,j}(u) = 
\left|
\begin{array}{cccc}
T_{1,1}(u) & \cdots &T_{1,i-1}(u)& T_{1,i}(u)\\
\vdots & \ddots &\vdots&\vdots\\
T_{i-1,1}(u) & \cdots & T_{i-1,i-1}(u)&T_{i-1,i}(u)\\
T_{j,1}(u) & \cdots & T_{j,i-1}(u)&
\hbox{\begin{tabular}{|c|}\hline$T_{j,i}(u)$\\\hline\end{tabular}}
\end{array}
\right|D_i(u)^{-1}.
\end{equation}
Since $E_{j-1}^{(1)} = T_{j-1,j}^{(1)}$ and $F_{j-1}^{(1)} = T_{j,j-1}^{(1)}$,
it follows easily that
\begin{equation}\label{eij}
E_{i,j}^{(r)} = 
[E_{i,j-1}^{(r)}, E_{j-1}^{(1)}],\qquad
F_{i,j}^{(r)} = 
[F_{j-1}^{(1)},F_{i,j-1}^{(r)}].
\end{equation}
for $i+1 < j \leq n$. 
Comparing (\ref{qd1})--(\ref{qd3}) 
with Lemma~\ref{qdet} we deduce:

\begin{Lemma}\label{reduce}
For all admissible $i$, we have that
\begin{itemize}
\item[(i)]
$D_i(u) = \psi_{i-1}(D_{1}(u)) = \psi_{i-1}(T_{1,1}(u))$;
\item[(ii)]
$E_i(u) = \psi_{i-1}(E_{1}(u)) = \psi_{i-1}(T_{1,1}(u)^{-1} T_{1,2}(u))$;
\item[(iii)]
$F_i(u) = \psi_{i-1}(F_{1}(u)) = \psi_{i-1}(T_{2,1}(u) T_{1,1}(u)^{-1})$.
\end{itemize}
\end{Lemma}

In particular Lemma~\ref{reduce}(i)
shows that the elements $D_i^{(r)}$ here are the same as
the elements denoted $D_{i;1,1}^{(r)}$ 
from $\S$\ref{slevi}, so they generate
the Cartan subalgebra $Y_{(1^n)}$.
Also let ${Y}_{(1^n)}^+$ resp. ${Y}_{(1^n)}^-$ denote the subalgebra
of ${Y}_n$ generated by the elements
$\{E_i^{(r)}\}_{i=1,\dots,n-1,r \geq 1}$
resp. $\{F_i^{(r)}\}_{i=1,\dots,n-1,r \geq 1}$.
In view of (\ref{eij}), 
all the elements $E_{i,j}^{(r)}$ belong
to $Y_{(1^n)}^+$ and all the elements $F_{i,j}^{(r)}$ belong to $Y_{(1^n)}^-$.
By applying the antiautomorphism $\tau$ to the factorization
(\ref{gfact}), one checks:
\begin{equation}\label{tau}
\tau(E_{i,j}(u))  = F_{i,j}(u),\qquad
\tau(F_{i,j}(u)) = E_{i,j}(u),\qquad
\tau(D_i(u)) = D_i(u).
\end{equation}
Hence, $\tau$ fixes $Y_{(1^n)}$ elementwise
and interchanges the subalgebras
$Y_{(1^n)}^+$ and $Y_{(1^n)}^-$.

Now we state the main theorem of the section.
This is essentially due to Drinfeld \cite{D3}; see the remark at the end
of the section for the precise relationship.

\begin{Theorem}\label{drinpres} The algebra $Y_n$ is generated by the elements
$\{D_i^{(r)}, \widetilde{D}_i^{(r)}\}_{1 \leq i \leq n, r \geq 0}$
and 
$\{E_i^{(r)}, F_i^{(r)}\}_{1 \leq i < n,r \geq 1}$
subject only to the following relations:
\begin{align}
D_i^{(0)} &= 1,\label{r0}\\
\sum_{t=0}^r D_i^{(t)} \widetilde D_i^{(r-t)} &= -\delta_{r,0},\label{r1}\\
[D_i^{(r)}, D_j^{(s)}] &=  0,\label{r2}\\
[E_i^{(r)},F_j^{(s)}] &= \delta_{i,j} 
\sum_{t=0}^{r+s-1} \widetilde D_{i}^{(t)}D_{i+1}^{(r+s-1-t)} ,\label{r3}\\
[D_i^{(r)}, E_j^{(s)}] &= (\delta_{i,j}-\delta_{i,j+1})
\sum_{t=0}^{r-1} D_i^{(t)} E_j^{(r+s-1-t)},\label{r4}\\
[D_i^{(r)}, F_j^{(s)}] &= (\delta_{i,j+1}-\delta_{i,j})
\sum_{t=0}^{r-1} F_j^{(r+s-1-t)}D_i^{(t)} ,\label{r5}
\end{align}\begin{align}
[E_i^{(r)}, E_i^{(s)}] &=
\sum_{t=1}^{s-1} E_i^{(t)} E_i^{(r+s-1-t)}
-\sum_{t=1}^{r-1} E_i^{(t)} E_i^{(r+s-1-t)},\label{r6}\\
[F_i^{(r)}, F_i^{(s)}] &=
\sum_{t=1}^{r-1} F_i^{(r+s-1-t)} F_i^{(t)}-\sum_{t=1}^{s-1} 
F_i^{(r+s-1-t)} F_i^{(t)},\label{r7}\\
 [E_i^{(r)}, E_{i+1}^{(s+1)}]&- [E_i^{(r+1)}, E_{i+1}^{(s)}] =
-E_i^{(r)} E_{i+1}^{(s)},\label{r8}\\
[F_i^{(r+1)}, F_{i+1}^{(s)}] &- [F_i^{(r)}, F_{i+1}^{(s+1)}] =
 -F_{i+1}^{(s)} F_i^{(r)},\label{r9}\\
[E_i^{(r)}, E_j^{(s)}] &= 0 \hspace{36.9mm}\text{ if }|i-j|> 1,\label{r10}\\
[F_i^{(r)}, F_j^{(s)}] &= 0 \hspace{36.9mm}\text{ if }|i-j|> 1,\label{r11}\\
[E_i^{(r)}, [E_i^{(s)}, E_j^{(t)}]] &+ 
[E_i^{(s)}, [E_i^{(r)}, E_j^{(t)}]] = 0 \quad\text{ if }|i-j|=1,\label{r12}\\
[F_i^{(r)}, [F_i^{(s)}, F_j^{(t)}]] &+ 
[F_i^{(s)}, [F_i^{(r)}, F_j^{(t)}]] = 0 \quad\text{ if }|i-j|=1,\label{r13}
\end{align}
for all admissible $i,j,r,s, t$.
\end{Theorem}

\begin{Remark}\rm
The relations (\ref{r6}) and (\ref{r7}) are equivalent to the relations
\begin{align}
[E_i^{(r)}, E_i^{(s+1)}] &- [E_i^{(r+1)}, E_i^{(s)}] =
E_i^{(r)} E_i^{(s)} + E_i^{(s)} E_i^{(r)},\label{r6b}\\
[F_i^{(r+1)}, F_i^{(s)}] &- [F_i^{(r)}, F_i^{(s+1)}] =
F_i^{(r)} F_i^{(s)} + F_i^{(s)} F_i^{(r)},\label{r7b}
\end{align}
respectively.
\end{Remark}

In the remainder of the section, we are going to write down a proof,
since we could not find one in the literature.
There are two parts to the proof:
first we must show that all these relations are satisfied in $Y_n$; 
second we must show that we have found enough relations.

Let us begin with some reductions to the first part of the proof.
We have already noted that 
the elements $\{D_i^{(r)}\}_{i=1,\dots,n, r \geq 1}$ commute, hence
the relations (\ref{r0})--(\ref{r2}) hold. 
Also by Lemma~\ref{reduce},
$D_i^{(r)} \in \psi_{i-1}(Y_1)$ and $E_j^{(s)} \in \psi_{j-1}(Y_2)$,
so Lemma~\ref{cent} implies that
(\ref{r4}) holds if either $i < j$ or $i > j+1$.
Similar reasoning shows that (\ref{r3}) holds
if $|i-j| > 1$ and (\ref{r10})--(\ref{r11}) hold always.
Having made these remarks, Lemma~\ref{reduce} and (\ref{tau})
reduces the verification of all the remaining
relations to checking the following special cases:
(\ref{r3}) with $i=1,j=1$ or $i=2,j=1$;
(\ref{r4}) with $i=1,j=1$ or $i=2,j=1$;
(\ref{r6}) with $i=1$; 
(\ref{r8}) with $i=1$;
(\ref{r12}) with $i=2,j=1$ or $i=1,j=2$.

\begin{Lemma}\label{goody2}
The following identities hold in $Y_2((u^{-1},v^{-1}))$:
\begin{itemize}
\item[(i)]
$(u-v) [D_1(u), E_1(v)] = D_1(u) (E_1(v)-E_1(u))$;
\item[(ii)]
$(u-v) [E_{1}(u), \widetilde{D}_{2}(v)]
 = (E_{1}(u) - E_{1}(v)) \widetilde D_{2}(v)$;
\item[(iii)]
$(u-v) [E_1(u), F_1(v)] = \widetilde D_1(v) D_2(v) - \widetilde D_1(u) D_2(u)$;
\item[(iv)]
$(u-v)[E_1(u), E_1(v)] = (E_1(v) - E_1(u))^2$.
\end{itemize}
\end{Lemma}

\begin{proof}
Compute the $e_{1,1} \otimes e_{1,2}$-, 
$e_{1,2} \otimes e_{2,2}$-, $e_{1,2} \otimes e_{2,1}$-
and $e_{1,2} \otimes e_{1,2}$-coefficients on each side of
(\ref{rmdefalt})
and rearrange the resulting four equations
to obtain the following identities:
\begin{itemize}
\item[(i)$'$]
$(u-v) [T_{1,1}(u), \widetilde T_{1,2}(v)] = T_{1,1}(u) \widetilde T_{1,2}(v)
+ T_{1,2}(u) \widetilde T_{2,2}(v)$;
\item[(ii)$'$]
$(u-v) [T_{1,2}(u), \widetilde T_{2,2}(v)] = 
T_{1,1}(u) \widetilde T_{1,2}(v)
+ T_{1,2}(u) \widetilde T_{2,2}(v)$;
\item[(iii)$'$]
$(u-v) [T_{1,2}(u), \widetilde T_{2,1}(v)] = T_{1,1}(u) \widetilde
T_{1,1}(v) + T_{1,2}(u) \widetilde T_{2,1}(v)
-\widetilde T_{2,1}(v) T_{1,2}(u)- \widetilde T_{2,2}(v) T_{2,2}(u)$;
\item[(iv)$'$]
$[T_{1,2}(u), \widetilde T_{1,2}(v)] = 0$.
\end{itemize}
We also note that
$$
\left(
\begin{array}{ll}
T_{1,1}(u) & T_{1,2}(u)\\
T_{2,1}(u) & T_{2,2}(u)
\end{array}
\right)
= 
\left(
\begin{array}{ll}
D_1(u)& D_1(u) E_1(u)\\
F_1(u) D_1(u)& D_2(u) + F_1(u) D_1(u) E_1(u)
\end{array}
\right)
$$
and that
$$
\left(
\begin{array}{cc}
\widetilde T_{1,1}(v) & \widetilde T_{1,2}(v)\\
\widetilde T_{2,1}(v) & \widetilde T_{2,2}(v)
\end{array}
\right)
= 
\left(
\begin{array}{cc}
\widetilde D_1(v) + E_1(v) \widetilde D_2(v) F_1(v)& - E_1(v) \widetilde D_2(v)\\
- \widetilde D_2(v) F_1(v)& \widetilde D_2(v)
\end{array}
\right).
$$
Substituting from these into (i)$'$ and 
using the known fact that $D_1(u)$ commutes with $\widetilde D_2(v)$
gives the identity
$$
(u-v) [D_1(u), E_1(v)] \widetilde D_2(v) = D_1(u) (E_1(v) - E_1(u)) \widetilde D_2(v).
$$
Multiplying on the right by $D_2(v)$ gives (i).
The deduction of (ii) from (ii)$'$ is entirely similar.

Next we deduce (iii) from (iii)$'$. 
Rewriting (i) and (ii) using $\tau$ gives that
\begin{align*}
E_1(v) \widetilde D_2(v)
+(u-v-1) E_1(u) \widetilde D_2(v)&=(u-v) \widetilde D_2(v) E_1(u)\\
F_1(u) D_1(u)
+(u-v-1)F_1(v) D_1(u) &= (u-v) D_1(u) F_1(v).
\end{align*}
Also rearranging (iii)$'$ gives
\begin{multline*}
D_1(u) \widetilde D_1(v)
+D_1(u) 
(E_1(v) \widetilde D_2(v)
+(u-v-1) E_1(u) \widetilde D_2(v) ) F_1(v)
=\\
D_2(u) \widetilde D_2(v)
+\widetilde D_2(v) 
(F_1(u) D_1(u)
+(u-v-1)F_1(v) D_1(u))E_1(u).
\end{multline*}
Now substituting the first two of these identities into the third
 and multiplying on the left by 
$\widetilde D_1(u) D_2(v)$ gives (iii).

Finally we must deduce (iv). By (iv)$'$, we have that
$$
D_1(u) E_1(u) E_1(v) \widetilde D_2(v) = E_1(v) D_1(u) \widetilde D_2(v) 
E_1(u).
$$
Multiply both sides by $(u-v)^2$ and use (i) and (ii) to move $D_1(u)$
to the left and $\widetilde D_2(v)$ to the right, then cancel the
leading $D_1(u)$'s and trailing $\widetilde D_2(v)$'s to get
$$
(u-v)^2 E_1(u) E_1(v) = 
\left((u-v) E_1(v) - E_1(v) + E_1(u)\right)
\left((u-v) E_1(u) + E_1(v) - E_1(u)\right).
$$
Hence,
\begin{itemize}
\item[(iv)$''$]
$(u-v)^2 [E_1(u), E_1(v)] = 
(E_1(v)-E_1(u))(E_1(u)-E_1(v))
+ (u-v) E_1(v) (E_1(v) - E_1(u))
+ (u-v) (E_1(u) - E_1(v)) E_1(u).$
\end{itemize}
Now 
subtract $(u-v) [E_1(u), E_1(v)]$ from both sides of (iv)$''$
to deduce that
$$
(u-v)(u-v-1) [E_1(u),E_1(v)]
=
(u-v-1) (E_1(v) - E_1(u))^2.
$$
Hence (iv) follows on dividing both sides by $(u-v-1)$.
\end{proof}

\begin{Lemma}\label{goody3}
The following identities hold in $Y_3((u^{-1}, v^{-1}))$:
\begin{itemize}
\item[(i)] $[E_1(u), F_2(v)] = 0$;
\item[(ii)] $(u-v) [E_1(u), E_2(v)] = E_1(u)E_2(v)-E_1(v)E_2(v) - E_{1,3}(u) + 
E_{1,3}(v)$;
\item[(iii)] $[E_{1,3}(u), E_2(v)] = E_2(v) [E_1(u), E_2(v)]$;
\item[(iv)] 
$[E_{1}(u), E_{1,3}(v) - E_{1}(v) E_{2}(v)] = -[E_{1}(u), E_{2}(v)] E_{1}(u)$.
\end{itemize}
\end{Lemma}

\begin{proof}
Arguing as in the proof of the previous lemma,
we compute the 
$e_{1,2} \otimes e_{3,2}$-,
$e_{1,2} \otimes e_{2,3}$-,
$e_{1,3} \otimes e_{2,3}$- and $e_{1,2} \otimes e_{1,3}$-coefficients of 
(\ref{rmdefalt}) respectively to obtain the identities
\begin{itemize}
\item[(i)$'$] $[D_1(u) E_1(u), \widetilde D_3(v) F_2(v)] = 0$
\item[(ii)$'$] $(u-v) [D_1(u) E_1(u), E_2(v) \widetilde D_3(v)] =
D_1(u) ( E_1(u)E_2(v) - E_1(v) E_2(v) + E_{1,3}(v) - E_{1,3}(u)) \widetilde D_3(v)$;
\item[(iii)$'$] $[D_1(u)E_{1,3}(u), E_{2}(v) \widetilde D_3(v)] = 0$;
\item[(iv)$'$] $[D_1(u) E_1(u), (E_{1,3}(v) - E_1(v) E_2(v)) \widetilde D_3(v)] = 0$.
\end{itemize}
Using commutation relations already derived, (i) and (ii)
follow easily from (i)$'$ and (ii)$'$. 
To prove (iii), we need one more identity.
By Lemma~\ref{goody2}(ii), we know that
$$
(u-v) [\widetilde D_3(v), E_2(u)] = (E_2(v) - E_2(u))\widetilde D_3(v).
$$
Considering $u^0$-coefficients gives
$[\widetilde D_3(v), E_2^{(1)}] = E_2(v) \widetilde D_3(v)$.
Hence, recalling (\ref{eij}),
\begin{align*}
[E_{1,3}(u), \widetilde D_3(v)] 
&=
[[E_{1}(u), E_{2}^{(1)}], \widetilde D_3(v)]
=
[E_1(u), [E_{2}^{(1)}, \widetilde D_3(v)]]\\
&=
-[E_1(u), E_2(v) \widetilde D_3(v)] = -[E_1(u), E_2(v)] \widetilde D_3(v).
\end{align*}
Now take (iii)$'$, cancel the leading $D_1(u)$ and then
simplify to get
\begin{align*}
[E_{1,3}(u), E_{2}(v) \widetilde D_3(v)]&=
[E_{1,3}(u), E_2(v)] \widetilde D_3(v)
+ E_2(v) [E_{1,3}(u), \widetilde D_3(v)]\\
&=
\left([E_{1,3}(u), E_2(v)] - E_2(v) [E_1(u), E_2(v)]  \right)
\widetilde D_3(v) = 0.
\end{align*}
This proves (iii).
For (iv), note to start with by considering the
$u^0$-coefficients of (ii) that
$[E_1^{(1)}, E_2(v)] = E_{1,3}(v) - E_1(v) E_2(v)$.
Lemma~\ref{goody2}(i) implies that
$[D_1(u), E_1^{(1)}] = D_1(u) E_1(u)$.
Now compute:
\begin{align*}
[D_1(u), E_{1,3}(v) - E_1(v) E_2(v)] &=
[D_1(u), [E_1^{(1)}, E_2(v)]]
=
[[D_1(u), E_1^{(1)}], E_2(v)] \\
&=
[D_1(u) E_1(u), E_2(v)] = D_1(u) [E_1(u), E_2(v)].
\end{align*}
Using this identity to rewrite (iv)$'$ we get
\begin{multline*}
[D_1(u) E_1(u), E_{1,3}(v) - E_1(v) E_2(v)] 
=
D_1(u) [E_1(u), E_2(v)] E_1(u)\\
+ D_1(u) [E_1(u), E_{1,3}(v) - E_1(v) E_2(v)]
= 0.
\end{multline*}
Now (iv) follows on cancelling $D_1(u)$.
\end{proof}

\begin{Lemma}\label{serre1}
The following relations hold:
\begin{itemize}
\item[(i)] $[[E_1(u), E_2(v)], E_2(v)] = 0$;
\item[(ii)] $[E_1(u), [E_1(u), E_2(v)]] = 0$.
\end{itemize}
\end{Lemma}

\begin{proof}
(i) Compute using Lemma~\ref{goody3}(ii) and (iii):
\begin{align*}
(u-v) 
[[E_1(u), E_2(v)], E_2(v)]&=
[E_1(u) E_2(v) - E_1(v) E_2(v) - E_{1,3}(u) + E_{1,3}(v), E_2(v)]\\
&=
[E_1(u),E_2(v)] E_2(v) 
- [E_1(v),E_2(v)] E_2(v)\\
&\qquad\qquad+ E_2(v) [E_1(v),E_2(v)] - E_2(v)[E_1(u),E_2(v)]\\
&=[[E_1(u),E_2(v)], E_2(v)] 
- [[E_1(v),E_2(v)], E_2(v)].
\end{align*}
We conclude that
$(u-v-1) [[E_1(u), E_2(v)], E_2(v)] = - [[E_1(v),E_2(v)], E_2(v)]$.
Now let $u = v+1$ to deduce that the right hand side equals zero, 
then divide by $(u-v-1)$ to proof the lemma.

(ii) Similar calculation using Lemma~\ref{goody3}(iv)
instead of (iii).
\end{proof}

\begin{Lemma}\label{serre2}
The following relations hold:
\begin{itemize}
\item[(i)] $[[E_1(u), E_2(v)], E_2(w)] + 
[[E_1(u), E_2(w)], E_2(v)] = 0$;
\item[(ii)] $[E_1(u), [E_1(v), E_2(w)]] + [E_1(v), [E_1(u), E_2(w)]] = 0$.
\end{itemize}
\end{Lemma}

\begin{proof}
(i)
We show that the expression 
$(u-v)(u-w)(v-w) [[E_1(u), E_2(v)], E_2(w)]$
is symmetric in $v$ and $w$. 
By Lemma~\ref{goody3}(ii) it equals
$$
(u-w)(v-w)
[E_1(u) E_2(v) - E_1(v) E_2(v) + E_{1,3}(v) - E_{1,3}(u), E_2(w)]
$$
Using Lemmas~\ref{goody3}(iii) and  \ref{serre1}(i) 
this equals
\begin{align*}
&(u-w)(v-w)[E_1(u), E_2(w)] E_2(v)
+(u-w)(v-w)E_1(u) [E_2(v), E_2(w)]\\
-\,&(u-w)(v-w)[E_1(v), E_2(w)]  E_2(v)
-(u-w)(v-w)E_1(v) [E_2(v), E_2(w)]\\
+\,&(u-w)(v-w) [E_1(v),E_2(w)]E_2(w) 
-(u-w)(v-w) [E_1(u), E_2(w)]E_2(w) 
\end{align*}
Now use Lemmas~\ref{goody3}(ii) and \ref{goody2}(iv) to expand the
commutators once more to get
\begin{align*}
&(v-w)(E_1(u) E_2(w)E_2(v) - E_1(w) E_2(w)E_2(v) + E_{1,3}(w)E_2(v) - E_{1,3}(u)  E_2(v))\\
+\,&(u-w) (E_1(u)E_2(v)^2 - E_1(u)E_2(v) E_2(w) - E_1(u)E_2(w) E_2(v) +E_1(u) E_2(w)^2)\\
-\,&(u-w)(E_1(v) E_2(w)E_2(v) - E_1(w)E_2(w)E_2(v) + E_{1,3}(w)E_2(v) - E_{1,3}(v)  E_2(v))\\
-\,&(u-w) (E_1(v)E_2(v)^2 - E_1(v)E_2(v) E_2(w) - E_1(v)E_2(w) E_2(v) + E_1(v)E_2(w)^2)\\
+\,&(u-w)  (E_1(v) E_2(w)E_2(w) - E_1(w) E_2(w)E_2(w) + E_{1,3}(w)E_2(w) - 
E_{1,3}(v)E_2(w))\\
-&(v-w) (E_1(u) E_2(w)E_2(w) - E_1(w) E_2(w)E_2(w) + E_{1,3}(w)E_2(w) - E_{1,3}(u)E_2(w)).
\end{align*}
Now open the parentheses and check that the resulting expression
is symmetric in $v$ and $w$ to complete the proof.

(ii) A similar calculation using Lemma~\ref{goody3}(iv)
instead of (iii) and Lemma~\ref{serre1}(ii) instead of (i)
shows that the expression 
$(u-v)(u-w)(v-w) [E_1(u), [E_1(v), E_2(w)]]$
is symmetric in $u$ and $v$. 
\end{proof}

Now we can verify the remaining relations needed for the first 
part of the proof. Note that 
\begin{equation}\label{todd}
(E_1(v) - E_1(u))/(u-v)
= \sum_{r,s \geq 1} E_1^{(r+s-1)} u^{-r} v^{-s}.
\end{equation}
Using this, divide both sides of the identity from Lemma~\ref{goody2}(i)
by $(u-v)$ and equate $u^{-r} v^{-s}$-coefficients on both sides to prove
(\ref{r4}) with $i=1,j=1$.
Next, multiplying
Lemma~\ref{goody2}(ii) on the left and right by $D_2(v)$
then swapping $u$ and $v$ gives the identity
\begin{equation}\label{newby}
(u-v) [D_2(u), E_1(v)] = -D_2(u)(E_1(v) - E_1(u)).
\end{equation}
Now argue using (\ref{todd}) 
again to deduce (\ref{r4}) with $i=2,j=1$ from this.
Similarly one gets (\ref{r3}) with $i=1,j=1$
from Lemma~\ref{goody2}(iii), 
(\ref{r3}) with $i=2,j=1$ from Lemma~\ref{goody3}(i),
(\ref{r6}) with $i=1$ from Lemma~\ref{goody2}(iv), 
(\ref{r8}) with $i=1$ from Lemma~\ref{goody3}(ii),
(\ref{r12}) with $i=2,j=1$ from Lemma~\ref{serre2}(i) 
and (\ref{r12}) with $i=1,j=2$ from Lemma~\ref{serre2}(ii).

\vspace{2mm}

Now we consider the second part of the proof.
Let $\widehat{Y}_n$ denote the algebra with generators and relations
as in the statement of Theorem~\ref{drinpres}.
For $1 \leq i < j \leq n$, define elements $E_{i,j}^{(r)},
F_{i,j}^{(r)} \in \widehat{Y}_n$
by the equations (\ref{eij}).
Let $\widehat{Y}_{(1^n)}$ resp. 
$\widehat{Y}_{(1^n)}^+$ resp. $\widehat{Y}_{(1^n)}^-$ denote the subalgebra
of $\widehat{Y}_n$ generated by the elements
$\{D_i^{(r)}\}_{i=1,\dots,n, r \geq 1}$ resp.
$\{E_i^{(r)}\}_{i=1,\dots,n-1,r \geq 1}$
resp. $\{F_i^{(r)}\}_{i=1,\dots,n-1,r \geq 1}$.
Define an ascending filtration
$$
\L_0 \widehat{Y}^+_{(1^n)} \subseteq \L_1
\widehat{Y}^+_{(1^n)} \subseteq \cdots
$$
on $\widehat{Y}^+_{(1^n)}$ by declaring that the generator $E_i^{(r)}$ 
is of degree $(r-1)$, i.e. $\L_d \widehat{Y}^+_{(1^n)}$ is the span of
all monomials in these generators of total degree at most $d$.
Let $\gr^\LL \widehat{Y}_{(1^n)}^+$ denote the associated graded algebra,
and let $e_{i,j;r} := \gr^\LL_{r} E_{i,j}^{(r+1)} 
\in \gr^\LL \widehat{Y}_{(1^n)}^+$
for each $1 \leq i < j \leq n$ and $r \geq 0$.

\begin{Lemma}\label{main2}
For $1 \leq i < j \leq n, 1 \leq \k < \l \leq n$
and $r,s \geq 0$,
we have that
$$
[e_{i,j;r}, e_{\k,\l;s}]
=  e_{i,\l;r+s}\delta_{\k,j} - 
\delta_{i,\l}e_{\k,j;r+s}.
$$
\end{Lemma}

\begin{proof}
By the defining relations for $\widehat{Y}_n$,
we have easily that
\begin{itemize}
\item[(i)] $[e_{i,i+1;r}, e_{j,j+1;s}] = 0$
if $|i-j| \neq 1$;
\item[(ii)] $[e_{i,i+1;r+1}, e_{j,j+1;s}] =
[e_{i,i+1;r}, e_{j,j+1;s+1}]$ if $|i-j| = 1$;
\item[(iii)] 
$[e_{i,i+1;r}, [e_{i,i+1;s}, e_{j,j+1;t}]] = 
-[e_{i,i+1;s}, [e_{i,i+1;r}, e_{j,j+1;t}]]$
if $|i-j|=1$.
\end{itemize}
We also have by definition 
that 
\begin{itemize}
\item[(iv)]
$e_{i,j;r} = 
[e_{i,j-1;r}, e_{j-1,j;0}]$
for $j > i+1$. 
\end{itemize}
Now we consider seven cases.
\begin{itemize}
\item[(1)] $j < \k$. Obviously, 
$[e_{i,j;r}, e_{\k,\l;s}] = 0$.
\item[(2)] $j=\k$. By (ii) and (iv),
$[e_{j-1,j;r}, e_{j,j+1;s}] = 
e_{j-1,j+1;r+s}$. 
Now bracket with
$e_{j+1,j+2;0}, \dots, e_{\l-1,\l;0}$
to deduce that
$[e_{j-1,j;r}, e_{j,\l;s}] = 
e_{j-1,\l;r+s}$.
Finally bracket with $e_{j-2,j-1;0},\dots,
e_{i,i+1;0}$.
\item[(3)] $i < \k, j=\l$.
Let us just show that $[e_{1,3;r}, e_{2,3;s}] = 0$, since the
general case is an easy consequence. 
Note that by (iii), $[e_{1,3;r}, e_{2,3;s}] =
[[e_{1,2;r}, e_{2,3;0}], e_{2,3;s}] 
= -[[e_{1,2;r}, e_{2,3;s}], e_{2,3;0}].$
By (ii) this equals
$-[[e_{1,2;r+s}, e_{2,3;0}], e_{2,3;0}]$ which is zero by (iii).
\item[(4)] $i = \k, j < \l$. Similar to (3).
\item[(5)] $i=\k,j=\l$. If $j = i+1$, we are done by (i);
else,
\begin{align*}
[e_{i,j;r}, e_{i,j;s}] &=
[[e_{i,j-1;r}, e_{j-1,j;0}], e_{i,j;s}]\\
&=
[[e_{i,j-1;r}, e_{i,j;s}], e_{j-1,j;0}]
+ [e_{i,j-1;r}, [e_{j-1,j;0}, e_{i,j;s}]]
\end{align*} 
which is
zero by (3) and (4).
\item[(6)] $i < \k < j < \l$.
We just show $[e_{1,3;r}, e_{2,4;s}] = 0$.
It equals
\begin{align*}
[[e_{1,2;r}, e_{2,3;0}],[e_{2,3;0}, e_{3,4;s}]] &=
[e_{2,3;0}, [[e_{1,2;r}, e_{2,3;0}],e_{3,4;s}]]\\&=
[e_{2,3;0}, [e_{1,2;r}, [e_{2,3;0},e_{3,4;s}]]]\\
&=
[[e_{2,3;0}, e_{1,2;r}],[e_{2,3;0}, e_{3,4;s}]]\\&=
-[[e_{1,2;r}, e_{2,3;0}],
[e_{2,3;0}, e_{3,4;s}]].
\end{align*}
Hence it is zero.
\item[(7)] $i < \k, \l < j$.
We just show $[e_{2,3;r}, e_{1,4;s}] = 0$.
It equals $[e_{2,3;r}, [e_{1,3;s}, e_{3,4;0}]] = [e_{1,3;s},  e_{2,4;r}]$, 
which is zero by (6). 
\end{itemize}
\end{proof}

\begin{Lemma} \label{span}
The algebra $\widehat{Y}_n$ 
is spanned by the set of monomials
in $\{D_i^{(r)}\}_{i=1,\dots,n,r \geq 1}
\cup\{E_{i,j}^{(r)}, F_{i,j}^{(r)}\}_{1 \leq i < j \leq n, r \geq 1}$, 
taken in some fixed order so that
$F$'s come before $D$'s and $D$'s come before $E$'s.
\end{Lemma}

\begin{proof}
Using Lemma~\ref{main2}, one shows easily that 
the associated graded algebra $\gr^\LL 
\widehat{Y}_{(1^n)}^+$ is spanned by the set of all ordered monomials
in the elements $\{e_{i,j;r}\}_{1 
\leq i < j \leq n, r \geq 0}$
taken in some fixed order. Hence $\widehat{Y}_{(1^n)}^+$ itself
is spanned by the corresponding monomials in 
$\{{E}_{i,j}^{(r)}\}_{1 \leq i < j \leq n, r \geq 1}$.
By the defining relations for $\widehat{Y}_n$, 
there is an antiautomorphism $\tau$
of $\widehat{Y}_n$ fixing each $D_i^{(r)}$ and interchanging
each $E_{i,j}^{(r)}$ with $F_{i,j}^{(r)}$.
It follows on applying $\tau$ that
the set of monomials
in $\{F_{i,j}^{(r)}\}_{1 \leq i < j \leq n}$
taken in some fixed order span $\widehat{Y}_{(1^n)}^-$, while
obviously the ordered monomials
in the elements $\{D_i^{(r)}\}_{i=1,\dots,n, r \geq 1}$ span
$\widehat{Y}_{(1^n)}$. Since by the defining relations
the natural multiplication map 
$\widehat{Y}_{(1^n)}^- \otimes \widehat{Y}_{(1^n)} \otimes \widehat{Y}_{(1^n)}^+ 
\rightarrow \widehat{Y}_n$ is surjective, the lemma follows.
\end{proof}

Now, the first part of the proof of Theorem~\ref{drinpres} above implies that
there is a surjective algebra homomorphism
$\theta:\widehat{Y}_n \rightarrow Y_n$ sending
$D_i^{(r)}, E_{i,j}^{(r)},F_{i,j}^{(r)}\in\widehat{Y}_n$ 
to the elements with the same name in $Y_n$.
To complete the proof of Theorem~\ref{drinpres}
we need to show that $\theta$ is an isomorphism.
This follows immediately from Lemma~\ref{ind} below, since it shows that
the images of the monomials that span $\widehat{Y}_n$
from Lemma~\ref{span} are linearly independent in $Y_n$ itself.

\begin{Lemma} \label{ind}
The set of monomials
in $\{D_i^{(r)}\}_{i=1,\dots,n,r \geq 1}
\cup\{E_{i,j}^{(r)}, F_{i,j}^{(r)}\}_{1 \leq i < j \leq n,r \geq 1}$
taken in some fixed order is linearly independent in $Y_n$.
\end{Lemma}

\begin{proof}
As explained at the end of $\S$\ref{srtt},
we can identify the associated graded algebra $\gr^\LL Y_n$ with
$U(\mathfrak{gl}_n[t])$, so that $\gr^\LL_{r} T_{i,j}^{(r+1)}$
is identified with $e_{i,j} t^{r}$.
It is easy to see from (\ref{qd1})--(\ref{qd3}) that under this identification
$\gr^\LL_{r} D_i^{(r+1)}$
resp. $\gr^\LL_{r} E_{i,j}^{(r+1)}$
resp. $\gr^\LL_{r} F_{i,j}^{(r+1)}$ is identified with
$e_{i,i} t^{r}$ resp. $e_{i,j} t^{r}$ resp. $e_{j,i} t^{r}$.
Hence by the PBW theorem for $U(\mathfrak{gl}_n[t])$,
the set of all monomials in
$$
\{\gr^\LL_{r} D_i^{(r+1)}\}_{i=1,\dots,n,r \geq 0}
\cup\{\gr^\LL_{r} E_{i,j}^{(r+1)},
\gr^\LL_{r} F_{i,j}^{(r+1)}\}_{1 \leq i < j \leq n, r \geq 0}
$$ 
taken in some fixed order forms a basis for $\gr^\LL Y_n$.
The lemma follows easily.
\end{proof}

This completes the proof of Theorem~\ref{drinpres}.
Let us also state the following theorem 
which was obtained in the course of the above proof; cf. \cite{Lev}.

\begin{Theorem}\label{triangular}
\begin{itemize}
\item[(i)] The set of all monomials in 
$\{D_i^{(r)}\}_{i=1,\dots,n, r \geq 1}$ taken in some fixed order
form a basis for $Y_{(1^n)}$.
\item[(ii)] The set of all monomials in 
$\{E_{i,j}^{(r)}\}_{1 \leq i < j \leq n, r \geq 1}$ taken in some fixed
order form a basis for $Y_{(1^n)}^+$.
\item[(iii)] The set of all monomials in 
$\{F_{i,j}^{(r)}\}_{1 \leq i < j \leq n, r \geq 1}$ taken in some fixed
order form a basis for $Y_{(1^n)}^-$.
\item[(iv)] The set of all monomials in 
$\{D_i^{(r)}\}_{i=1,\dots,n, r \geq 1}\cup\{E_{i,j}^{(r)},
F_{i,j}^{(r)}\}_{1 \leq i < j \leq n, r \geq 1}$ taken in some fixed
order form a basis for $Y_n$.
\end{itemize}
\end{Theorem}

\begin{Remark}\label{drinrem}\rm
Let us explain the relationship between 
the presentation given in Theorem~\ref{drinpres} and 
Drinfeld's presentation from \cite{D3}, since there are
some additional shifts in $u$.
Actually, the latter is a presentation for the subalgebra
$$
Y(\mathfrak{sl}_n) = 
\{x \in Y_n\:|\:\mu_f(x) = x \hbox{ for all }
f(u) \in 1 + \C[[u^{-1}]]\};
$$
see \cite[Definition 2.14]{MNO}.
Define $\kappa_{i,k}, \xi_{i,k}^{\pm}$ for $i=1,\dots,n-1$
and $k \geq 0$ from the equations
 \begin{align}
\kappa_i(u) = \sum_{k \geq 0} \kappa_{i,k}u^{-k-1}
&:=
\textstyle
1+\widetilde{D}_i\left(u -\frac{i-1}{2}\right) 
D_{i+1}\left(u-\frac{i-1}{2}\right),\\
\xi_i^+(u) = \sum_{k \geq 0} \xi_{i,k}^+ u^{-k-1}
&:= 
\textstyle
E_{i}\left(u-\frac{i-1}{2}\right),\\
\xi_i^-(u) = \sum_{k \geq 0} \xi_{i,k}^- u^{-k-1}&:=
\textstyle
F_{i}\left(u-\frac{i-1}{2}\right).
\end{align}
One can check by equating coefficients in the identities
from Lemmas~\ref{goody2}, \ref{goody3} and \ref{serre2} that 
these elements generate $Y(\mathfrak{sl}_n)$ subject to 
the Drinfeld relations, namely:
\begin{align}\label{dr1}
[\kappa_{i,k}, \kappa_{j,l}] &= 0, \\
[\xi^+_{i,k}, \xi^-_{j,l}] &= \delta_{i,j} \kappa_{i,k+l},\\
[\kappa_{i,0}, \xi^{\pm}_{j,l}] &= \pm a_{i,j} \xi^{\pm}_{j,l},\\
[\kappa_{i,k},\xi^\pm_{j,l+1}] - [\kappa_{i,k+1}, \xi^{\pm}_{j,l}]
&= \pm \frac{a_{i,j}}{2} (\kappa_{i,k} \xi^{\pm}_{j,l} + \xi^{\pm}_{j,l} \kappa_{i,k}),\\
[\xi^\pm_{i,k}, \xi^\pm_{j,l+1}] - [\xi^\pm_{i,k+1}, \xi^\pm_{j,l}]
&= \pm \frac{a_{i,j}}{2} (\xi^\pm_{i,k} \xi^{\pm}_{j,l} + \xi^\pm_{j,l} \xi^\pm_{i,k}),\label{eg}\\
i \neq j, N = 1 - a_{i,j} \Rightarrow
\operatorname{Sym} &[\xi_{i,k_1}^{\pm}, [\xi_{i,k_2}^{\pm}, \cdots
[\xi_{i,k_N}^{\pm}, \xi_{j,l}^\pm]\cdots]] = 0
\label{drn}
\end{align}
where $(a_{i,j})_{1 \leq i,j < n}$ denotes the Cartan matrix of type
$A_{n-1}$ indexed in the standard way and $\operatorname{Sym}$ denotes
symmetrization with respect to $k_1,\dots,k_N$.
For example, let us verify (\ref{eg}) in the case $j=i+1$ and the sign is $+$:
applying Lemma~\ref{goody3}(ii)
\begin{align*}
(u-v) [\xi_i^+(u), \xi_{i+1}^+(v)]
&=
\left(\left(u-{\textstyle\frac{i-1}{2}}\right) - \left(v-{\textstyle\frac{i}{2}}\right) - {\textstyle\frac{1}{2}}\right)
\left[E_i\left(u-{\textstyle\frac{i-1}{2}}\right), E_{i+1}\left(v-{\textstyle\frac{i}{2}}\right)\right]\\
&=
E_i\left(u-{\textstyle\frac{i-1}{2}}\right) E_{i+1}\left(v-{\textstyle\frac{i}{2}}\right)
- E_i\left(v-{\textstyle\frac{i}{2}}\right) E_{i+1}\left(v-{\textstyle\frac{i}{2}}\right)\\
&\phantom{=}\qquad- E_{i,i+2}\left(u-{\textstyle\frac{i-1}{2}}\right) + E_{i,i+2} \left(v-{\textstyle\frac{i}{2}}\right)\\
&\phantom{=}\qquad- {\textstyle\frac{1}{2}} \left[E_i\left(u-\textstyle{\frac{i-1}{2}}\right), E_{i+1}\left(v-{\textstyle\frac{i}{2}}\right)\right]\\
&=
{\textstyle\frac{1}{2}}\left(\xi_i^+(u) \xi_{i+1}^+(v)+\xi_{i+1}^+(v)\xi_i^+(u)\right)
- E_i\left(v-{\textstyle\frac{i}{2}}\right) E_{i+1}\left(v-{\textstyle\frac{i}{2}}\right)\\
&\phantom{=}\qquad- E_{i,i+2}\left(u-{\textstyle\frac{i-1}{2}}\right) + E_{i,i+2} \left(v-{\textstyle\frac{i}{2}}\right).
\end{align*}
Now equate $u^{-k-1}v^{-l-1}$-coefficients on both sides to get
$$
[\xi_{i,k+1}^+, \xi_{i+1,l}^+] - [\xi_{i,k}^+, \xi_{i+1,l+1}^+]
= {\textstyle\frac{1}{2}}\left(
\xi_{i,k}^+ \xi_{i+1,l}^+ + \xi_{i+1,l}^+ \xi_{i,k}^+
\right)
$$
as required.
({\em Beware}: the relations (\ref{dr1})--(\ref{drn}) are actually not exactly the same
as the relations in \cite{D3} --- 
one needs to swap $\xi_{i,k}^+$ and $\xi_{i,k}^-$
and replace $\kappa_{i,k}$ with $-\kappa_{i,k}$ to get those.
The reason for the difference is that we have chosen to work with
the opposite presentation to Drinfeld throughout.)
\end{Remark}

\section{Parabolic subalgebras}\label{sparabolic}

Now we are ready to prove Theorems A and B stated in the introduction.
Fix throughout the section a tuple $\nu = (\nu_1,\dots,\nu_m)$
of non-negative integers summing to $n$. 
Note there is going to be 
some overlap between the notation here and that of the
previous section, which is the special case $\nu = (1^n)$ of the 
present definitions. When necessary, we will add an additional
superscript $\nu$ to our notation to avoid any ambiguity as $\nu$ varies.
Factor the $n \times n$ matrix $T(u)$ as 
\begin{equation}\label{gfact2}
T(u) = F(u) D(u) E(u)
\end{equation}
for unique {\em block matrices}
$$
D(u) = \left(
\begin{array}{cccc}
D_{1}(u) & 0&\cdots&0\\
0 & D_{2}(u) &\cdots&0\\
\vdots&\vdots&\ddots&\vdots\\
0&0 &\cdots&D_{m}(u)
\end{array}
\right),
$$
$$
E(u) = 
\left(
\begin{array}{cccc}
I_{\nu_1} & E_{1,2}(u) &\cdots&E_{1,m}(u)\\
0 & I_{\nu_2} &\cdots&E_{2,m}(u)\\
\vdots&\vdots&\ddots&\vdots\\
0&0 &\cdots&I_{\nu_m}
\end{array}
\right),\:
F(u) = \left(
\begin{array}{cccc}
I_{\nu_1} & 0 &\cdots&0\\
F_{1,2}(u) & I_{\nu_2} &\cdots&0\\
\vdots&\vdots&\ddots&\vdots\\
F_{1,m}(u)&F_{2,m}(u) &\cdots&I_{\nu_m}
\end{array}
\right),
$$
where $D_a(u) = (D_{a;i,j}(u))_{1 \leq i,j \leq \nu_a}$, 
$E_{a,b}(u)= (E_{a,b;i,j}(u))_{1 \leq i \leq \nu_a, 1 \leq j \leq \nu_b}$ 
and $F_{a,b}(u) = (F_{a,b;i,j}(u))_{1 \leq i \leq \nu_b, 1 \leq j \leq 
\nu_a}$ are  $\nu_a \times \nu_a$, 
$\nu_a \times \nu_b$
and  $\nu_b \times\nu_a$ matrices, respectively.
Also define the $\nu_a \times \nu_a$ matrix 
$\widetilde{D}_a(u) = (\widetilde{D}_{a;i,j}(u))_{1 \leq i,j \leq \nu_a}$
by $\widetilde{D}_a(u) := -D_a(u)^{-1}$.
The entries of these matrices define power series
$D_{a;i,j}(u) = \sum_{r \geq 0} D_{a;i,j}^{(r)} u^{-r}$, 
$\widetilde{D}_{a;i,j}(u) = \sum_{r \geq 0}\widetilde{D}_{a;i,j}^{(r)} 
u^{-r}$, 
$E_{a,b;i,j}(u) = \sum_{r \geq 1} E_{a,b;i,j}^{(r)} u^{-r}$ and 
$F_{a,b;i,j}(u) = \sum_{r \geq 1} F_{a,b;i,j}^{(r)} u^{-r}$.
We let 
$E_{a;i,j}(u) = \sum_{r \geq 1} E_{a;i,j}^{(r)} u^{-r} 
:= E_{a,a+1;i,j}(u)$ and 
$F_{a;i,j}(u) = \sum_{r \geq 1} F_{a;i,j}^{(r)} u^{-r} 
:= F_{a,a+1;i,j}(u)$ for short.

Like before, there are explicit descriptions of all these
elements in terms of quasi-determinants. To write them down,
write the matrix $T(u)$ in block form as
$$
T(u) = \left(
\begin{array}{lll}
{^\nu}T_{1,1}(u)&\cdots&{^\nu}T_{1,m}(u)\\
\vdots&\ddots&\cdots\\
{^\nu}T_{m,1}(u)&\cdots&{^\nu}T_{m,m}(u)\\
\end{array}
\right)
$$
where ${^\nu}T_{a,b}(u)$ is a $\nu_a \times \nu_b$ matrix.
Then, recalling the notation (\ref{qdetdef}),
\begin{equation}
\!\!\!\!\!\!
\!\!\!\!\!\!
\!\!\!\!\!\!
\!\!\!D_a(u) = 
\left|
\begin{array}{cccc}
{^\nu}T_{1,1}(u) & \cdots & {^\nu}T_{1,a-1}(u)&{^\nu}T_{1,a}(u)\\
\vdots & \ddots &\vdots&\vdots\\
{^\nu}T_{a-1,1}(u)&\cdots&{^\nu}T_{a-1,a-1}(u)&{^\nu}T_{a-1,a}(u)\\
{^\nu}T_{a,1}(u) & \cdots & {^\nu}T_{a,a-1}(u)&
\hbox{\begin{tabular}{|c|}\hline${^\nu}T_{a,a}(u)$\\\hline\end{tabular}}
\end{array}
\right|,
\end{equation}
\begin{equation}
E_{a,b}(u) = 
D_a(u)^{-1} 
\left|\begin{array}{cccc}
{^\nu}T_{1,1}(u) & \cdots &{^\nu}T_{1,a-1}(u)& {^\nu}T_{1,b}(u)\\
\vdots & \ddots &\vdots&\vdots\\
{^\nu}T_{a-1,1}(u) & \cdots & {^\nu}T_{a-1,a-1}(u)&{^\nu}T_{a-1,b}(u)\\
{^\nu}T_{a,1}(u) & \cdots & {^\nu}T_{a,a-1}(u)&
\hbox{\begin{tabular}{|c|}\hline${^\nu}T_{a,b}(u)$\\\hline\end{tabular}}
\end{array}
\right|,
\end{equation}
\begin{equation}
F_{a,b}(u) = 
\left|
\begin{array}{cccc}
{^\nu}T_{1,1}(u) & \cdots &{^\nu}T_{1,a-1}(u)& {^\nu}T_{1,a}(u)\\
\vdots & \ddots &\vdots&\vdots\\
{^\nu}T_{a-1,1}(u) & \cdots & {^\nu}T_{a-1,a-1}(u)&{^\nu}T_{a-1,a}(u)\\
{^\nu}T_{b,1}(u) & \cdots & {^\nu}T_{b,a-1}(u)&
\hbox{\begin{tabular}{|c|}\hline${^\nu}T_{b,a}(u)$\\\hline\end{tabular}}
\end{array}
\right|D_a(u)^{-1}.
\end{equation}
It follows in particular from these descriptions that
for $b > a+1$ and $1 \leq i \leq \nu_a, 1 \leq j \leq \nu_b$,
\begin{equation}
E_{a,b;i,j}^{(r)} = [E_{a,b-1;i,k}^{(r)}, E_{b-1;k,j}^{(1)}],
\qquad
F_{a,b;j,i}^{(r)} = [F_{b-1;j,k}^{(1)}, F_{a,b-1;k,i}^{(r)}],\label{ter}
\end{equation}
for any $1 \leq k \leq \nu_{b-1}$.
We also get the analogue of Lemma~\ref{reduce}:

\begin{Lemma}\label{reduce2}
Fix $a \geq 1$ and let $\bar\nu := (\nu_a,\nu_{a+1},\dots,\nu_m)$.
Then, for all admissible $i,j$,
\begin{itemize}
\item[(i)]
${^\nu}D_{a;i,j}(u) = \psi_{\nu_1+\cdots+\nu_{a-1}}({^{\bar\nu}}D_{1;i,j}(u))$;
\item[(ii)]
${^\nu}E_{a;i,j}(u) = 
\psi_{\nu_1+\cdots+\nu_{a-1}}({^{\bar\nu}} E_{1;i,j}(u))$;
\item[(iii)]
${^\nu} F_{a;i,j}(u) = 
\psi_{\nu_1+\cdots+\nu_{a-1}}({^{\bar\nu}}F_{1;i,j}(u))$.
\end{itemize}
\end{Lemma}

In particular Lemma~\ref{reduce2}(i)
shows that the elements $D_{a;i,j}^{(r)}$ here are the same as
the generators 
of the standard Levi subalgebra $Y_\nu$
introduced at the end of $\S$\ref{slevi}, so they satisfy
the relations (\ref{levirel}).
We also let
$Y_{\nu}^+$ resp. $Y_{\nu}^-$ denote the subalgebra generated by
$\{E_{a;i,j}^{(r)}\}_{1\leq a < m, 1 \leq i \leq \nu_a, 1 \leq 
j \leq \nu_{a+1}, r \geq 1}$
resp. $\{F_{a;i,j}^{(r)}\}_{1\leq a < m, 1 \leq i \leq \nu_{a+1}, 1 \leq 
j \leq \nu_{a}, r \geq 1}$.
The antiautomorphism $\tau$ has the properties
\begin{align}
\tau(D_{a;i,j}(u)) &= D_{a;j,i}(u),\label{tau1}\\
\tau(E_{a,b;i,j}(u)) &= F_{a,b;j,i}(u),\\
\tau(F_{a,b;i,j}(u)) &= E_{a,b;j,i}(u),\label{tau3}
\end{align}
so it leaves $Y_\nu$ invariant and 
interchanges $Y_{\nu}^+$ and $Y_{\nu}^-$.
Finally, as in the introduction, set $Y_\nu^\sharp := Y_\nu Y_\nu^+$
and $Y_\nu^\flat := Y_\nu^- Y_\nu$, giving the 
standard parabolic subalgebras of $Y_n$ of shape $\nu$.

\begin{Remark}\label{new}\rm
Here is an alternative definition of the
algebras $Y_\nu. Y_\nu^\sharp$ and $Y_\nu^\flat$ just in terms of the
Drinfeld generators from Theorem~\ref{drinpres}.
First, $Y_\nu$ is the subalgebra of $Y_n$ generated
by all $\{D_i^{(r)}\}_{1 \leq i \leq n,r \geq 1}$
together with $\{E_i^{(r)}, F_i^{(r)}\}_{r \geq 1}$
for $i \in \{1,\dots,n\} - \{\nu_1,\nu_1+\nu_2,\dots,\nu_1+\cdots+\nu_m\}$.
Then $Y_\nu^\sharp$ resp. $Y_\nu^\flat$ is the subalgebra generated
by $Y_\nu$ together with all remaining
elements of 
$\{E_i^{(r)}\}_{1 \leq i < n,r \geq 1}$
resp.
$\{F_i^{(r)}\}_{1 \leq i < n,r \geq 1}$.
\end{Remark}

We have now defined the elements $D_{a;i,j}^{(r)},
E_{a;i,j}^{(r)}$ and $F_{a;i,j}^{(r)}$ appearing in Theorems A and B
stated in the introduction.
We are ready to explain the proofs of these theorems.
Actually, the argument runs
almost exactly parallel to the proofs of Theorems~\ref{drinpres} and 
\ref{triangular} given in the previous section. 
As before, there are two parts: first, to show all the
relations (\ref{pr1})--(\ref{pr14}) from Theorem A hold; 
second, to show we have enough relations by constructing the PBW
bases described in Theorem B.

For the first part, one uses Lemma~\ref{reduce2}, (\ref{levirel}) 
and (\ref{tau1})--(\ref{tau3}) to reduce the problem
to checking the following special cases:
(\ref{pr6}) with $a=1,b=1$ or $a=2,b=1$;
(\ref{pr4}) with $a=1,b=1$ or $a=2,b=1$;
(\ref{pr7}) with $a=1$; 
(\ref{pr9}) with $a=1$;
(\ref{pr11}) with $a=1,b=2$;
(\ref{pr13}) with $a=2,b=1$ or $a=1,b=2$.
These special cases may be deduced from the following four lemmas
by equating coefficients. Note these lemmas 
are the exact analogues of Lemmas~\ref{goody2}--\ref{serre2}.

\begin{Lemma}\label{pgoody2}
Suppose $m=2$, i.e. $\nu = (\nu_1,\nu_2)$.
The following identities hold for all admissible $h,i,j,k$:
\begin{itemize}
\item[(i)]
$(u-v) [D_{1;i,j}(u), E_{1;\k,\l}(v)] = 
D_{1;i,p}(u) (E_{1;p,\l}(v) - E_{1;p,\l}(u))
\delta_{\k,j}$;
\item[(ii)]
$(u-v) [E_{1;i,j}(u), \widetilde{D}_{2;\k,\l}(v)]
 = (E_{1;i,q}(u) - E_{1;i,q}(v)) \widetilde D_{2;q,\l}(v)\delta_{\k,j}$;
\item[(iii)]
$(u-v) [E_{1;i,j}(u), F_{1;\k,\l}(v)] = 
\widetilde D_{1;i,\l}(v) D_{2;\k,j}(v) - \widetilde D_{1;i,\l}(u) D_{2;\k,j}(u)$;
\item[(iv)]
$(u-v)[E_{1;i,j}(u), E_{1;\k,\l}(v)] = (E_{1;i,\l}(v) - E_{1;i,\l}(u))
(E_{1;\k,j}(v) - E_{1;\k,j}(u))$.
\end{itemize}
(Here, $p$ resp. $q$ should be summed over $1,\dots,\nu_1$ 
resp. $1,\dots,\nu_2$.)
\end{Lemma}

\begin{proof}
Compute the $e_{i,j} \otimes e_{\k,\nu_1+\l}$-,
$e_{i,\nu_1+j} \otimes e_{\nu_1+\k,\nu_1+\l}$-,
$e_{i,\nu_1+j} \otimes e_{\nu_1+\k,\l}$-
and $e_{i,\nu_1+j} \otimes e_{\k,\nu_1+\l}$-coefficients on each side of
(\ref{rmdefalt})  and then rearrange the resulting identities
like we did in the proof of Lemma~\ref{goody2} above to
obtain:
\begin{itemize}
\item[(i)$'$] $(u-v) [D_{1;i,j}(u), E_{1;\k,q}(v) \widetilde{D}_{2;q,\l}(v)]
= 
\\\hspace{20mm}
D_{1;i,p}(u) (E_{1;p,q}(v) - E_{1;p,q}(u)) \widetilde{D}_{2;q,\l}(v) 
\delta_{\k,j}$;
\item[(ii)$'$] 
$(u-v) [D_{1;i,p}(u) E_{1;p,j}(u), \widetilde{D}_{2;\k,\l}(v)]
=
\\\hspace{20mm}
D_{1;i,p}(u)(E_{1;p,q}(u) - E_{1;p,q}(v))\widetilde{D}_{2;q,\l}(v)\delta_{\k,j};$
\item[(iii)$'$] $
(u-v)[D_{1;i,p}(u) E_{1;p,j}(u), \widetilde{D}_{2;\k,q}(v) F_{1;q,\l}(v)]
=\\
\hspace{20mm}
\delta_{i,\l} \widetilde{D}_{2;\k,q}(v) 
\big\{
D_{2;q,j}(u) + (F_{1;q,p}(u) - F_{1;q,p}(v)) D_{1;p,p'}(u) E_{1;p',j}(u)
\big\}\\\qquad\qquad
-
D_{1;i,p}(u) \big\{ 
\widetilde{D}_{1;p,\l}(v) + (E_{1;p,q}(v) - E_{1;p,q}(u)) \widetilde{D}_{2;q,q'}(v) F_{1;q',\l}(v)
\big\} \delta_{\k,j};$
\item[(iv)$'$] $
[D_{1;i,p}(u) E_{1;p,j}(u), E_{1;\k,q}(v) \widetilde{D}_{2;q,\l}(v)] =
0$.
\end{itemize}
(Here, $p,p'$ resp. $q,q'$ should also be summed over $1,\dots,\nu_1$
resp. $1,\dots,\nu_2$.)
Now (i), (ii) and (iii) are deduced from (i)$'$, (ii)$'$ and (iii)$'$
by simplifying exactly like we did in the proof of Lemma~\ref{goody2}.
It turns out to be more difficult than before to 
deduce (iv) from (iv)$'$ so we explain this part of the argument
more carefully.
As before, one rewrites (iv)$'$ using (i) and (ii) to obtain:
\begin{itemize}
\item[(iv)$''$]
$(u-v)^2 [E_{1;i,j}(u), E_{1;\k,\l}(v)]
=
(E_{1;i,j}(v) - E_{1;i,j}(u)) (E_{1;\k,\l}(u) - E_{1;\k,\l}(v))
+(u-v) E_{1;\k,j}(v) (E_{1;i,\l}(v) - E_{1;i,\l}(u))
+(u-v) (E_{1;i,\l}(u) - E_{1;i,\l}(v)) E_{1;\k,j}(u).
$
\end{itemize}
Now we deduce (iv) from this.
For a power series $X$ in $Y_n[[u^{-1},v^{-1}]]$, let us
write $\{X\}_d$ for the homogeneous component of $X$ of
total degree $d$ in the variables $u^{-1}$ and $v^{-1}$.
We show by induction on $d=1,2,\dots$ that
$$
(u-v) \{[E_{1;i,j}(u), E_{1;\k,\l}(v)]\}_{d+1}
=
\{(E_{1;i,\l}(v) - E_{1;i,\l}(u))(E_{1;\k,j}(v) - E_{1;\k,j}(u))\}_d.
$$
For the base case $d=1$, applying $\{.\}_0$ to 
(iv)$''$ shows
$(u-v)^2 \{[E_{1;i,j}(u), E_{1;\k,\l}(v)]\}_{2} = 0$,
hence 
$(u-v) \{[E_{1;i,j}(u), E_{1;\k,\l}(v)]\}_{2} = 0$ as required.
For the induction step, assume the statement is true for $d > 1$.
Apply $\{.\}_{d}$ to (iv)$''$
to get that
\begin{multline*}
(u-v)^2 \{[E_{1;i,j}(u), E_{1;\k,\l}(v)]\}_{d+2}
=
(u-v) 
\{E_{1;\k,j}(v) (E_{1;i,\l}(v) - E_{1;i,\l}(u))\}_{d+1}\\
\qquad\qquad\qquad\qquad\qquad\qquad\:\:- (u-v) \{(E_{1;i,\l}(v) - E_{1;i,\l}(u)) E_{1;\k,j}(u)\}_{d+1}\\
- \{(E_{1;i,j}(v) - E_{1;i,j}(u)) (E_{1;\k,\l}(v) - E_{1;\k,\l}(u))\}_d.
\end{multline*}
Now use the induction hypothesis, together with the identity
$\{[E_{1;\k,j}(v), E_{1;i,\l}(v)]\}_{d+1} =0$ which follows
by dividing both sides of the induction hypothesis 
by $(u-v)$ then setting $u=v$, to rewrite the right hand side
to deduce that
\begin{multline*}
(u-v)^2 \{[E_{1;i,j}(u), E_{1;\k,\l}(v)]\}_{d+2}
=
(u-v) 
\{(E_{1;i,\l}(v) - E_{1;i,\l}(u))E_{1;\k,j}(v) 
\}_{d+1}\\
- (u-v) \{(E_{1;i,\l}(v) - E_{1;i,\l}(u)) E_{1;\k,j}(u)\}_{d+1}.
\end{multline*}
Dividing both sides by $(u-v)$ completes the proof of the induction step.
\end{proof}

\begin{Lemma}\label{pgoody3}
Suppose $m=3$, i.e. $\nu = (\nu_1,\nu_2,\nu_3)$.
The following identities hold for all admissible $g,h,i,j,k$:
\begin{itemize}
\item[(i)] $[E_{1;i,j}(u), F_{2;\k,\l}(v)] = 0$;
\item[(ii)] $(u-v) [E_{1;i,j}(u), E_{2;\k,\l}(v)] = 
(E_{1;i,q}(u) E_{2;q,\l}(v)-E_{1;i,q}(v) E_{2;q,\l}(v)
- E_{1,3;i,\l}(u) + E_{1,3;i,\l}(v)) \delta_{\k,j}$;
\item[(iii)] 
$[E_{1,3;i,j}(u), E_{2;\k,\l}(v)] = E_{2;\k,j}(v) [E_{1;i,g}(u), E_{2;g,\l}(v)]$;
\item[(iv)] 
$[E_{1;i,j}(u), E_{1,3;\k,\l}(v) - E_{1;\k,q}(v) E_{2;q,\l}(v)]
= -[E_{1;i,g}(u), E_{2;g,\l}(v)] E_{1;\k,j}(u)$.
\end{itemize}
(Here, $q$ should
be summed over $1,\dots,\nu_2$.)
\end{Lemma}

\begin{proof}
One computes the $e_{i,\nu_1+j} \otimes e_{\nu_1+\nu_2+\k,\nu_1+\l}$-,
$e_{i,\nu_1+j} \otimes e_{\nu_1+\k,\nu_1+\nu_2+\l}$-,
$e_{i,\nu_1+\nu_2+j} \otimes e_{\nu_1+\k,\nu_1+\nu_2+\l}$- and 
$e_{i,\nu_1+j} \otimes e_{\k,\nu_1+\nu_2+\l}$-coefficients of 
(\ref{rmdefalt}) respectively like in the proof of Lemma~\ref{goody3}
to obtain the identities
\begin{itemize}
\item[(i)$'$] $[D_{1;i,p}(u) E_{1;p,j}(u), \widetilde D_{3;\k,r}(v) 
F_{2;r,\l}(v)] = 0$;
\item[(ii)$'$] $(u-v)
[D_{1;i,p}(u) E_{1;p,j}(u), E_{2;\k,r}(v) \widetilde D_{3;r,\l}(v)] =
D_{1;i,p}(u) ( E_{1;p,q}(u)E_{2;q,r}(v) - \\
E_{1;p,q}(v) E_{2;q,r}(v) + E_{1,3;p,r}(v) - E_{1,3;p,r}(u)) \widetilde D_{3;r,\l}(v) \delta_{\k,j}$;
\item[(iii)$'$] $[D_{1;i,p}(u)E_{1,3;p,j}(u), E_{2;\k,r}(v) \widetilde D_{3;r,\l}(v)] = 0$;
\item[(iv)$'$] $[D_{1;i,p}(u) E_{1;p,j}(u), (E_{1,3;\k,r}(v) - E_{1;\k,q}(v) E_{2;q,r}(v)) \widetilde D_{3;r,\l}(v)] = 0$.
\end{itemize}
(Here, $p,q$ and $r$ sum over $1,\dots,\nu_1,1,\dots,\nu_2$ and $1,\dots,\nu_3$
respectively.) Now (i)--(iv) are deduced from (i)$'$--(iv)$'$ by copying
the arguments from the proof of Lemma~\ref{goody3}.
\end{proof}

\begin{Lemma}\label{pserre1}
Suppose $m=3$, i.e. $\nu = (\nu_1,\nu_2,\nu_3)$.
The following identities hold for all admissible $f,g,h,i,j,k$:
\begin{itemize}
\item[(i)] $[[E_{1;i,j}(u), E_{2;\k,\l}(v)], E_{2;f,g}(v)] = 0$;
\item[(ii)] $[E_{1;i,j}(u), [E_{1;\k,\l}(u), E_{2;f,g}(v)]] = 0$.
\end{itemize}
\end{Lemma}

\begin{proof}
Dividing both sides of Lemma~\ref{pgoody2}(iv) by $(u-v)$ then
setting $v=u$ shows that $[E_{a;i,j}(u), E_{a;\k,\l}(u)] = 0$.
Given this and Lemma~\ref{pgoody3}(ii), (i) is obvious
unless $f=\k=j$ and (ii) is obvious unless
$f=\l=j$. Now the proof in these cases is completed exactly
like the proof of Lemma~\ref{serre1}.
\end{proof}

\begin{Lemma}\label{pserre2}
Suppose $m=3$, i.e. $\nu = (\nu_1,\nu_2,\nu_3)$.
The following identities hold for all admissible $f,g,h,i,j,k$:
\begin{itemize}
\item[(i)] $[[E_{1;i,j}(u), E_{2;\k,\l}(v)], E_{2;f,g}(w)] + 
[[E_{1;i,j}(u), E_{2;\k,\l}(w)], E_{2;f,g}(v)] = 0$;
\item[(ii)] $[E_{1;i,j}(u), [E_{1;\k,\l}(v), E_{2;f,g}(w)]] + [E_{1;i,j}(v), [E_{1;\k,\l}(u), E_{2;f,g}(w)]] = 0$.
\end{itemize}
\end{Lemma}

\begin{proof}
Show that 
$(u-v)(u-w)(v-w) [[E_{1;i,j}(u), E_{2;j,\l}(v)], E_{2;f,g}(w)]$
is symmetric in $v$ and $w$ 
and that
$(u-v)(u-w)(v-w) [E_{1;i,j}(u), [E_{1;\k,\l}(v), E_{2;\l,g}(w)]]$
is symmetric in $u$ and $v$,
following the argument of Lemma~\ref{serre2} exactly.
\end{proof}

Now we consider the second part of the proof.
Let $\widehat{Y}_n$ denote the algebra with generators and relations
as in the statement of Theorem~A.
Define elements $E_{a,b;i,j}^{(r)},
F_{a,b;j,i}^{(r)} \in \widehat{Y}_n$
by the equations (\ref{ter}).
We need to check that these definitions are independent of the particular 
choice of 
$k$. 
Well, given $1 \leq k,k' \leq \nu_{b-1}$ with $k \neq k'$, we have that
$[E_{a,b-1;i,k}^{(r)}, E_{b-1;k',j}^{(s)}]
=
0$ by (\ref{pr11}).
Bracketing with $D_{b-1;k,k'}^{(1)}$ and using (\ref{pr4}), one deduces that
\begin{equation}\label{indofk}
[E_{a,b-1;i,k}^{(r)}, E_{b-1;k,j}^{(s)}] =
[E_{a,b-1;i,k'}^{(r)}, E_{b-1;k',j}^{(s)}]
\end{equation}
as required to verify that the definition of the
elements $E_{a,b;i,j}^{(r)}$ is independent of the choice of $k$.
A similar argument shows that the definition of the elements
$F_{a,b;j,i}^{(r)}$ is independent of $k$ too.

Let $\widehat{Y}_\nu$,
$\widehat{Y}_\nu^+$ and $\widehat{Y}_{\nu}^-$ denote the subalgebras
of $\widehat{Y}_n$ generated by the 
$D$'s, $E$'s and $F$'s respectively.
By the first part of the proof, there is a surjective
homomorphism 
$\theta:\widehat{Y}_n \rightarrow Y_n$
sending $\widehat{Y}_\nu$ onto $Y_\nu$ and
$\widehat{Y}_\nu^{\pm}$ onto $Y_{\nu}^{\pm}$.
We just need to show that $\theta$ is an isomorphism. 
This is done just like in the previous section by exhibiting
a set of monomials that span $\widehat{Y}_n$ whose image in $Y_n$
is linearly independent. We just explain the key step, namely,
the analogue of Lemma~\ref{main2} allowing one to construct
the spanning set for $\widehat{Y}_{\nu}^+$. Given this, the rest of our
earlier argument extends without further
complication to complete the proof.
Define a filtration
$$
\L_0 \widehat{Y}^+_\nu \subseteq \L_1 \widehat{Y}^+_\nu \subseteq \cdots
$$
of $\widehat{Y}^+_\nu$ by declaring that the generators 
$E_{a;i,j}^{(r)}$ are of degree $(r-1)$. 
Let $\gr^\LL \widehat{Y}_\nu^+$ denote the associated graded algebra.
Letting $n_a := \nu_1+\cdots+\nu_{a-1}$ for short,
define $$
e_{n_a+i,n_b+j;r} := \gr^\LL_{r} 
E_{a,b;i,j}^{(r+1)} 
\in \gr^\LL \widehat{Y}_\nu^+
$$
for each $1 \leq a < b \leq m, 1 \leq i \leq \nu_a, 1 \leq j \leq \nu_b$
and $r\geq 0$. Then:

\begin{Lemma}\label{pmain2}
For $1 \leq a < b \leq m, 1 \leq c < d \leq m$, $r,s \geq 0$ 
and all admissible $h,i,j,k$,
we have that
$$
[e_{n_a+i,n_b+j;r}, 
e_{n_c+\k,n_d+\l;s}]
=  e_{n_a+i,n_d+\l;r+s}\delta_{n_c+\k,n_b+j} - 
\delta_{n_a+i,n_d+\l}e_{n_c+\k,n_b+j;r+s}.
$$
\end{Lemma}

\begin{proof}
Like in the proof of Lemma~\ref{main2}, we split into seven
cases: (1) $b < c$; (2) $b=c$; (3) $a < c, b = d$;
(4) $a = c, b < d$; (5) $a = c, b = d$;
(6) $a < c < b < d$; (7) $a < c; d < b$.
Since the analysis of each of the cases is very similar to 
Lemma~\ref{main2}, we just illustrate the idea
with the two hardest situations, both of which 
require the Serre relations.

First we check for case (3) that 
$[e_{n_1+i,n_3+j;r},e_{n_2+\k,n_3+\l;s}] = 0$.
For any $1 \leq g \leq \nu_2$,
we have by (\ref{ter}) and the images of the relations
(\ref{pr9}) and (\ref{pr13})
in $\gr^\LL \widehat{Y}_\nu$ that
\begin{align*}
[e_{n_1+i,n_3+j;r},e_{n_2+\k,n_3+\l;s}] &=
[[e_{n_1+i,n_2+g;r}, e_{n_2+g,n_3+j;0}],
e_{n_2+\k,n_3+\l;s}]\\&=
-[[e_{n_1+i,n_2+g;r}, e_{n_2+g,n_3+j;s}],
e_{n_2+\k,n_3+\l;0}]\\
&=-[[e_{n_1+i,n_2+g;r+s}, e_{n_2+g,n_3+j;0}],
e_{n_2+\k,n_3+\l;0}] = 0.
\end{align*}

Second we check for case (6) that $[e_{n_1+i,n_3+j;r},
e_{n_2+\k,n_4+\l;s}] = 0$.
By the case (1), (\ref{ter}), (\ref{pr9}) and (\ref{pr13}), we have that
\begin{align*}
[e_{n_1+i,n_3+j;r},e_{n_2+\k,n_4+\l;s}] &=
[[e_{n_1+i,n_2+\k;r}, e_{n_2+\k,n_3+j;0}],
[e_{n_2+\k,n_3+j;0}, e_{n_3+j,n_4+\l;s}]]
\\&
=[e_{n_2+\k,n_3+j;0}, [[e_{n_1+i,n_2+\k;r}, e_{n_2+\k,n_3+j;0}],e_{n_3+j,n_4+\l;s}]]\\
&=[e_{n_2+\k,n_3+j;0}, [e_{n_1+i,n_2+\k;r}, [e_{n_2+\k,n_3+j;0},
e_{n_3+j,n_4+\l;s}]]]\\&=[[e_{n_2+\k,n_3+j;0}, e_{n_1+i,n_2+\k;r}], [e_{n_2+\k,n_3+j;0},e_{n_3+j,n_4+\l;s}]]\\
&=-[[e_{n_1+i,n_2+\k;r}, e_{n_2+\k,n_3+j;0}], [e_{n_2+\k,n_3+j;0},e_{n_3+j,n_4+\l;s}]]\\&=-[e_{n_1+i,n_3+j;r},e_{n_2+\k,n_4+\l;s}]
\end{align*}
Hence it is zero.
\end{proof}

This completes the proof of Theorems A and B.

\section{Centers and centralizers}\label{scent}

In this section, we compute the centralizer in $Y_n$ of the standard
Levi subalgebra $Y_\nu$. 
The argument depends on the following auxiliary lemma, which is a 
generalization 
of \cite[Proposition 2.12]{MNO}; the proof given here is based on
the argument in {\em loc. cit.}.

\begin{Lemma}\label{aux}
Let $\mathfrak h$ be a reductive 
subalgebra of a finite dimensional Lie algebra $\mathfrak g$ over $\C$.
Let $\mathfrak c$ be the centralizer of $\mathfrak h$ in $\mathfrak g$.
Then, the centralizer of $U(\mathfrak h[t])$
in $U(\mathfrak g[t])$ is equal to $U(\mathfrak{c}[t])$.
\end{Lemma}

\noindent
(We believe that 
the word ``reductive'' is unnecessary here, but we did not
find a proof without it.)

\begin{proof}
The symmetrization map
$S(\mathfrak g[t]) \rightarrow U(\mathfrak g[t])$ is an isomorphism of
$\mathfrak h[t]$-modules. Using this, it suffices to show that the 
space of invariants of $\mathfrak h[t]$ acting on
$S(\mathfrak g[t])$ is $S(\mathfrak c[t])$.
Since $\mathfrak h$ is reductive, we can pick an 
$\ad \mathfrak h$-stable complement $\mathfrak c'$
to $\mathfrak c$ in $\mathfrak g$. 
Let $x_1,\dots,x_m$ be a basis for $\mathfrak c'$ and let
$x_{m+1},\dots,x_n$ be a basis for $\mathfrak c$.
Let $z$ be an $\mathfrak h[t]$-invariant in $S(\mathfrak g[t])$.
Define
$h \geq 0$ to be minimal such that 
$z$ has the form
$$
z = \sum_{d} z_d (x_1 t^h)^{d_1} \cdots (x_m t^h)^{d_m}
$$
summing over
$d = (d_1,\dots,d_m)$ with $d_1,\dots,d_m \geq 0$,
where the coefficients $z_d$ are polynomials in the variables
$x_i t^k$ for $i=1,\dots,m$ and $0 \leq k < h$ 
together with the variables $x_i t^k$ for $i=m+1,\dots,n$
and $k\geq 0$.
Pick a basis $y_1,\dots,y_r$ for $\mathfrak h$ and 
let
$[y_i, x_j] = \sum_{k=1}^m c_{i,j,k} x_k$
for each $j=1,\dots,m$.
Acting on $z$ with $y_i t \in \mathfrak h[t]$ and taking the coefficient
of $x_k t^{h+1}$ gives the equation
$$
\sum_d z_d \sum_{j=1}^m c_{i,j,k}
d_j (x_1 t^h)^{d_1} \cdots (x_j t^h)^{d_j-1}
\cdots (x_m t^h)^{d_m} = 0,
$$
for each $i=1,\dots,r$ and $k=1,\dots,m$.
Now fix $d = (d_1,\dots,d_m)$ with $d_1,\dots,d_m \geq 0$.
Taking the coefficient of $(x_1 t^h)^{d_1} \cdots (x_m t^h)^{d_m}$
in our equation gives 
\begin{align*}\label{have}
\sum_{j=1}^m
c_{i,j,k} (d_j+1)
z_{d+\delta_j} = 0\qquad&(i=1,\dots,r, k=1,\dots,m)\\\intertext{where 
$d+\delta_j$ denotes the tuple $(d_1,\dots,d_j+1,\dots,d_m)$.
Since 
$\mathfrak h$ has no non-trivial invariants in $\mathfrak c'$,
the system of linear equations
$[y_i, \sum_{j=1}^m \lambda_j x_j] = 0 \:\:(i=1,\dots,r)$
has only the trivial solution $\lambda_1=\cdots=\lambda_m =0$. 
Equivalently, the system of equations}
\sum_{j=1}^m c_{i,j,k} \lambda_j = 0 \qquad&(i=1,\dots,r, k = 1,\dots,m)
\end{align*}
has only the trivial solution $\lambda_1=\cdots=\lambda_m = 0$ too.
We deduce that $(d_j+1)z_{d+\delta_j} = 0$ for each 
$j=1,\dots,m$. Hence $z_d = 0$ for all non-zero $d$,
which implies by the minimality of the choice of $h$ 
that $h = 0$ hence that $z \in S(\mathfrak c[t])$.
\end{proof}

Now, working once more in terms of the usual Drinfeld generators
from $\S$\ref{sdrinfeld}, define
\begin{equation}\label{cdef}
C_n(u) = \sum_{r \geq 0} C_n^{(r)} u^{-r} :=
D_1(u)D_2(u-1) \cdots D_n(u-n+1).
\end{equation}
The importance of the elements $C_n^{(r)}$ is due to the following 
theorem; cf. \cite[Theorem 2.13]{MNO}.
We remark that this theorem implies in particular 
that the commutative subalgebra $Y_{(1^n)}$ of $Y_n$ is generated
by the centers $Z(Y_1), Z(Y_2), \dots, Z(Y_n)$ of the nested
subalgebras $Y_1 \subset Y_2 \subset\cdots\subset Y_n$ from (\ref{tower}).

\begin{Theorem}\label{shortly}
The elements $C_n^{(1)},C_n^{(2)},\dots$
are algebraically independent and generate the center $Z(Y_n)$.
\end{Theorem}

\begin{proof}
First we check that the $C_n^{(r)}$ are central.
For this, it suffices to 
show that $[D_i(u) D_{i+1}(u-1), E_i(v)] = 0 =
[D_i(u) D_{i+1}(u-1), F_i(v)]$ for each $i=1,\dots,n-1$.
Actually we just need to check the first equality, 
since the second then follows on applying $\tau$.
By Lemma~\ref{goody2}(i),
\begin{align*}
(u-v) E_i(v) D_i(u) &= (u-v-1) D_i(u) E_i(v)+D_i(u) E_i(u)\\\intertext{By (\ref{newby}),}
(u-v-1) E_i(v) D_{i+1}(u-1) &=(u-v) D_{i+1}(u-1)E_i(v) - D_{i+1}(u-1)E_i(u-1).\\\intertext{Hence, setting $v=u$,}
E_i(u) D_{i+1}(u-1) &= D_{i+1}(u-1)E_i(u-1).
\end{align*}
Now calculate 
$(u-v) E_i(v) D_i(u) D_{i+1}(u-1)$ using these identities to 
show that it equals
$(u-v) D_i(u) D_{i+1}(u-1)E_i(v)$.
Hence $[D_i(u) D_{i+1}(u-1), E_i(v)] = 0$.

Now we complete the proof by following
the argument of \cite[Theorem 2.13]{MNO}.
Recall the filtration (\ref{filt2}) of $Y_n$,
with associated graded algebra $\gr^\LL Y_n = U(\mathfrak{gl}_n[t])$.
Let
$z = e_{1,1}+\cdots+e_{n,n} \in \mathfrak{gl}_n$.
One checks from the definition
(\ref{cdef}) that
$$
\gr_{r-1}^\LL C_n^{(r)} = z t^{r-1}.
$$
By Lemma~\ref{aux} (taking $\mathfrak h = \mathfrak g = \mathfrak{gl}_n$)
the center of $U(\mathfrak{gl}_n[t])$ is
freely generated by the elements 
$\{zt^r\:|\:r \geq 0\}$.
The theorem now follows on combining these two observations with the fact already proved that
each $C_n^{(r)}$ belongs to $Z(Y_n)$.
\end{proof}

We now use essentially the same argument to prove
the following theorem, which is a well known variation on a result of Olshanskii
\cite[$\S$2.1]{Ol}.

\begin{Theorem}
The centralizer of $Y_m$ in $Y_{m+n}$
is equal to $Z(Y_m)\psi_m(Y_n)$.
\end{Theorem}

\begin{proof}
Lemma~\ref{cent} shows that
$Z(Y_m) \psi_m(Y_n)$ centralizes $Y_m$,
so we just need to show that the centralizer is no larger.
Consider the associated graded algebra $\gr^\LL Y_{m+n} 
= U(\mathfrak{gl}_{m+n}[t])$. Since
$$
\gr_{r-1}^\LL \psi_m(T_{i,j}^{(r)}) = 
e_{m+i,m+j} t^{r-1},
$$
we have that $\gr^\LL Y_m = U(\mathfrak{gl}_m[t])$ (where 
$\mathfrak{gl}_m$ 
is embedded into the top left corner of $\mathfrak{gl}_{m+n}$)
and $\gr^\LL \psi_m(Y_n) = U(\mathfrak{gl}_n[t])$ 
(where $\mathfrak{gl}_n$ is
embedded into the bottom right corner of $\mathfrak{gl}_{m+n}$).
By Lemma~\ref{aux}, the centralizer of $U(\mathfrak{gl}_m[t])$ in
$U(\mathfrak{gl}_{m+n}[t])$ is equal to $Z(U(\mathfrak{gl}_m[t])) 
U(\mathfrak{gl}_n[t])$. The theorem follows.
\end{proof}

\begin{Corollary}
Let $\nu = (\nu_1,\dots,\nu_m)$ be a tuple of non-negative integers summing
to $n$. The centralizer of the Levi subalgebra $Y_\nu$ in $Y_n$ is equal
to $Z(Y_\nu)$.
\end{Corollary}

\begin{proof}
Proceed by induction on $m$, the case $m=1$ being vacuous.
By (\ref{comp}) $Y_\nu = Y_{\nu_1} \psi_{\nu_1}(Y_{\bar\nu})$
where $\bar\nu = (\nu_2,\dots,\nu_m)$.
By the theorem, the centralizer in $Y_n$ of $Y_{\nu_1}$
is $Z(Y_{\nu_1}) \psi_{\nu_1}(Y_{\nu_2+\cdots+\nu_m})$.
By induction, the centralizer in
$Y_{\nu_2+\cdots+\nu_m}$ of $Y_{\bar\nu}$ is $Z(Y_{\bar\nu})$.
Hence, the centralizer of $Y_\nu$ in $Y_n$
is $Z(Y_{\nu_1}) \psi_{\nu_1}(Z(Y_{\bar\nu})) = Z(Y_\nu)$.
(Alternatively one can prove the corollary directly
using Lemma~\ref{aux} once more.)
\end{proof}

\begin{Corollary}
$Y_{(1^n)}$ is a maximal commutative subalgebra of $Y_n$.
\end{Corollary}

\begin{proof}
By the previous corollary, $Y_{(1^n)}$ 
is its own centralizer.
\end{proof}

\section{Quantum determinants}\label{squantum}

In the literature, Drinfeld generators are usually expressed
in terms {quantum determinants}, rather than the quasi-determinants
used up to now. In this section we complete the picture by
relating quasi-determinants to quantum determinants.
We begin by introducing quantum determinants following \cite[$\S$2]{MNO}.
Fix $d \geq 1$ and let $A_d \in M_n^{\otimes d}$
denote the antisymmetrization operator, 
i.e. the endomorphism
$$
v_1 \otimes \cdots \otimes v_d \mapsto \sum_{\pi \in S_d}
\sgn(\pi) v_{\pi 1} \otimes \cdots \otimes v_{\pi d}
$$
of the natural space $(\C^n)^{\otimes d}$ that
$M_n^{\otimes d}$ acts on.
Note that $A_d^2 = (d!)A_d$.
We have the following fundamental identity
\begin{multline}\label{fundid}
A_d^{[1,\dots,d]} T^{[1,d+1]}(u)T^{[2,d+1]}(u-1) \cdots
T^{[d,d+1]}(u-d+1)
=\\ 
T^{[d,d+1]}(u-d+1)\cdots 
T^{[2,d+1]}(u-1) T^{[1,d+1]}(u) 
A_d^{[1,\dots,d]} 
\end{multline}
equality written in $M_n^{\otimes d} \otimes Y_n[[u^{-1}]]$;
see \cite[Proposition 2.4]{MNO}.
For tuples  $\bi = (i_1,\dots,i_d)$ and $\bj = (j_1,\dots,j_d)$
of integers from $\{1,\dots,n\}$, 
the {\em quantum determinant}
$T_{\bi, \bj}(u) \in Y_n[[u^{-1}]]$ is defined to be the
coefficient of
$e_{\bi, \bj} = e_{i_1, j_1} \otimes \cdots\otimes e_{i_d,j_d} \in 
\Mat_n^{\otimes d}$
on either side of the equation (\ref{fundid}).
Explicit computation using the left and the right hand sides 
of (\ref{fundid}) respectively
gives that
\begin{align}\label{lfull}
T_{\bi,\bj}(u) &= \sum_{\pi \in S_d} \sgn(\pi)
T_{i_{\pi 1},j_{1}}(u)
T_{i_{\pi 2},j_{2}}(u-1)
\cdots T_{i_{\pi d}, j_{d}}(u-d+1)\\\label{rfull0}
&= 
\sum_{\pi \in S_d} \sgn(\pi)
T_{i_d,j_{\pi d}}(u-d+1) 
\cdots 
T_{i_{2}, j_{\pi 2}}(u-1)
T_{i_1, j_{\pi 1}}(u),
\end{align}
where $S_d$ is the symmetric group.
It is obvious from these formulae that
\begin{equation}\label{perm}
T_{\bi \cdot \pi, \bj}(u) =
\sgn(\pi) T_{\bi, \bj}(u) =T_{\bi, \bj \cdot \pi}(u)
\end{equation}
for any permutation $\pi \in S_d$ (acting naturally on the tuples
$\bi,\bj$ by place permutation).
Using (\ref{perm}) one obtains further variations
on the formulae (\ref{lfull})--(\ref{rfull0}) as in \cite[Remark 2.8]{MNO},
for instance:
\begin{equation}\label{rfull}
T_{\bi,\bj}(u) =
\sum_{\pi \in S_d} \sgn(\pi)
T_{i_1,j_{\pi 1}}(u-d+1) 
T_{i_2,j_{\pi 2}}(u-d+2) 
\cdots 
T_{i_d, j_{\pi d}}(u).
\end{equation}
The following properties of quantum determinants are easily
derived from (\ref{lfull})--(\ref{perm}).
\begin{align}\label{tauprop2}
\tau(T_{\bi,\bj}(u)) &= T_{\bj, \bi}(u),\\
\sigma(T_{\bi,\bj}(u)) &= T_{\bi,\bj}(-u+d-1).
\label{Sprop}
\end{align}
In the special case $\bi=\bj = (1,\dots,n)$, we denote the
quantum determinant $T_{\bi,\bj}(u)$ instead by $C_n(u)$, i.e.
\begin{equation}\label{newcdef}
C_n(u) := T_{(1,\dots,n), (1,\dots,n)}(u).
\end{equation}
We will show in Theorem~\ref{cid} below that this agrees with the definition
(\ref{cdef}), hence the coefficients of the series $C_n(u)$ generate
the center of $Y_n$, but we do not know this yet.

The next few results taken from 
\cite{NT} describe the effect of the maps
$\Delta$ and $S$ on quantum determinants.
Actually we do not need the first of these here, but include it for 
the sake of completeness.

\begin{Lemma}
Let $\bi,\bj$ be $d$-tuples of distinct integers from $\{1,\dots,n\}$.
Then,
$$
\Delta(T_{\bi,\bj}(u)) = \sum_{\bk} T_{\bi,\bk}(u) \otimes T_{\bk,\bj}(u)
$$
where the sum is over all $\bk = (k_1,\dots,k_d)$ 
with $1 \leq k_1 < \cdots < k_d \leq n$.
\end{Lemma}

\begin{proof}
See \cite[Proposition 1.11]{NT}.
\end{proof}

\begin{Lemma}\label{sl}
Let $\bi,\bj$ be $d$-tuples of distinct integers from $\{1,\dots,n\}$.
Choose $\bi' = (i_{d+1},\dots,i_n)$ and $\bj' = (j_{d+1},\dots,j_n)$
so that $\{i_1,\dots,i_n\} = 
\{j_1,\dots,j_{n}\} = \{1,\dots,n\}$,
and let $\eps$ denote the sign of the permutation
$(i_1,\dots,i_n) \mapsto (j_1,\dots,j_n)$. Then,
$$
\omega(T_{\bi,\bj}(u)) =
\eps  C_n(-u+n-1)^{-1} T_{\bj', \bi'}(-u+n-1).
$$
\end{Lemma}

\begin{proof}
This is proved in \cite[Lemma 1.5]{NT} but for the opposite algebra,
so we repeat the argument once more.
By the identity (\ref{fundid}) and the definition (\ref{newcdef}),
we have that
$$
A_n^{[1,\dots,n]} T^{[1,n+1]}(u) T^{[2,n+1]}(u-1) \cdots T^{[n,n+1]}(u-n+1)
=
C_n(u)^{[n+1]} A_n^{[1,\dots,n]};
$$
see also \cite[Proposition 2.5]{MNO}.
Hence,
\begin{multline*}
A_n^{[1,\dots,n]} T^{[1,n+1]}(u)\cdots 
T^{[d,n+1]}(u-d+1)
=\\
(-1)^{n-d}C_n(u)^{[n+1]} A_n^{[1,\dots,n]}
\widetilde{T}^{[n,n+1]}(u-n+1)
\cdots
\widetilde{T}^{[d+1,n+1]}(u-d).
\end{multline*}
Now equate the $e_{(i_1,\dots,i_d,i_n,\dots,i_{d+1}),
(j_1,\dots,j_d,i_n,\dots,i_{d+1})}$-coefficients on each side and
use (\ref{lfull}) to deduce that
$$
T_{\bi,\bj}(u) = \eps C_n(u) \omega(T_{\bj',\bi'}(-u+n-1)).
$$
The lemma follows on making some obvious substitutions.
\end{proof}

\begin{Corollary}\label{S1}
In the notation of Lemma~\ref{sl},
$S(T_{\bi,\bj}(u)) = \eps C_n(u+n-d)^{-1} T_{\bj', \bi'}(u+n-d)$.
\end{Corollary}

\begin{proof}
Recalling that $S = \omega\circ\sigma$,
this is a consequence of (\ref{Sprop}) and Lemma~\ref{sl}.
\end{proof}

\begin{Corollary}\label{s2}
$S^2(T_{i,j}(u)) = C_n(u+n)^{-1}T_{i,j}(u+n) C_n(u+n-1)$.
\end{Corollary}

\begin{proof}
Apply Corollary~\ref{S1} twice.
\end{proof}

Now we describe the embedding $\psi_m:Y_n \hookrightarrow Y_{m+n}$
from $\S$\ref{slevi} in terms of quantum determinants.

\begin{Lemma}\label{gr}
Let $\bi, \bj$
be $d$-tuples of distinct integers from $\{1,\dots,n\}$.
Then,
$$
\psi_m(T_{\bi,\bj}(u)) = C_m(u+m)^{-1}T_{m\#\bi,m\#\bj}(u+m)  
$$
where $m\#\bi$ 
denotes the $(m+d)$-tuple  $(1,\dots,m,m+i_1,\dots,m+i_d)$
and $m\#\bj$ is defined similarly.
\end{Lemma}

\begin{proof}
Calculate using (\ref{pdidef}) and Lemma~\ref{sl} twice.
\end{proof}

Now we can verify that $C_n(u)$
as defined in this section is the same as the
earlier definition (\ref{cdef}).
Note this theorem is the Yangian analogue of
\cite[Theorem 7.24]{GKLLRT}.

\begin{Theorem}\label{cid}
$C_n(u) = D_1(u) D_2(u-1) \cdots D_n(u-n+1)$.
\end{Theorem}

\begin{proof}
Recalling that $D_i(u) = \psi_{i-1}(T_{1,1}(u))$, 
Lemma~\ref{gr} implies that
$D_i(u) = C_{i-1}(u+i-1)^{-1} C_i(u+i-1)$.
The lemma follows easily from this by induction.
\end{proof}

Finally, we apply Lemma~\ref{gr} once more to
express the Drinfeld generators
from $\S$\ref{sdrinfeld}
in terms of quantum determinants; cf. \cite[Theorem B.15]{Iohara}
or \cite[Proposition 3.2]{Crampe}.

\begin{Theorem} \label{newd} For $i \geq 1$,
\begin{itemize}
\item[(i)] $D_i(u) = 
T_{(1,\dots,i-1), (1,\dots,i-1)}(u+i-1)^{-1}
T_{(1,\dots,i),(1,\dots,i)}(u+i-1)$;
\item[(ii)] $E_{i}(u) = T_{(1,\dots,i), (1,\dots,i)}(u+i-1)^{-1}
T_{(1,\dots,i), (1,\dots,i-1,i+1)}(u+i-1)$;
\item[(iii)] $F_{i}(u) = 
T_{(1,\dots,i-1,i+1), (1,\dots,i)}(u+i-1)T_{(1,\dots,i), 
(1,\dots,i)}(u+i-1)^{-1}$.
\end{itemize}
\end{Theorem}

\begin{proof}
Calculate using Lemmas~\ref{reduce} and \ref{gr} and the 
formulae (\ref{tau}) and (\ref{tauprop2}).
\end{proof}

\begin{Remark}\rm
Using Theorem~\ref{newd}, we can also express the
generating functions $\kappa_i(u)$ and $\xi_i^{\pm}(u)$
from Remark~\ref{drinrem} in terms of quantum determinants:
\begin{align*}
\kappa_i(u) &= \textstyle 1 - a_i(u)^{-1}a_i(u+1)^{-1} a_{i-1}\left(u + \frac{1}{2}\right)
a_{i+1}\left(u+\frac{1}{2}\right),\\
\xi_i^+(u) &= a_i(u)^{-1} b_i(u),\\
\xi_i^-(u) &= c_i(u)a_i(u)^{-1},
\end{align*}
where
$a_i(u) = T_{(1,\dots,i),(1,\dots,i)}\left(u+\frac{i-1}{2}\right)$,
$b_i(u) = T_{(1,\dots,i),(1,\dots,i-1,i+1)}\left(u+\frac{i-1}{2}\right)$ and
$c_i(u) = T_{(1,\dots,i-1,i+1),(1,\dots,i)}\left(u+\frac{i-1}{2}\right)$.
Allowing for the fact that we are working with the opposite presentation to
Drinfeld's (cf. Remark~\ref{drinrem}) this is the ``simpler isomorphism'' recorded immediately before
\cite[Theorem 2]{D3}.
\end{Remark}

\ifx\undefined\bysame
\newcommand{\bysame}{\leavevmode\hbox to3em{\hrulefill}\,}
\fi

\end{document}